\documentclass[mathpazo]{cicp}

\usepackage{graphicx}
\usepackage{mathtools}
\usepackage{subcaption}
\usepackage{color}

\newtheorem{defn}{Definition}

\begin{document}
%%%%% title : short title may not be used but TITLE is required.
% \title{TITLE}
% \title[short title]{TITLE}
\title{Splitting physics-informed neural networks for inferring the dynamics of integer- and fractional-order neuron models}

%%%%% author(s) :
% single author:
% \author[name in running head]{AUTHOR\corrauth}
% [name in running head] is NOT OPTIONAL, it is a MUST.
% Use \corrauth to indicate the corresponding author.
% Use \email to provide email address of author.
% \footnote and \thanks are not used in the heading section.
% Another acknowlegments/support of grants, state in Acknowledgments section
% \section*{Acknowledgments}
% \author[O.~Author]{Only Author\corrauth}
% \address{School of Mathematical Sciences, Beijing Normal University,
% Beijing 100875, P.R. China}
% \email{{\tt author@email} (O.~Author)}

% multiple authors:
% Note the use of \affil and \affilnum to link names and addresses.
% The author for correspondence is marked by \corrauth.
% use \emails to provide email addresses of authors
% e.g. below example has 3 authors, first author is also the corresponding
%      author, author 1 and 3 having the same address.
\author[Zhang Z R et.~al.]{Simin Shekarpaz\affil{1},
      Fanhai Zeng\affil{2}, and George Karniadakis\affil{1}\comma\corrauth}
\address{\affilnum{1}\ Division of Applied Mathematics, Brown University, Providence, RI 02912, USA. \\
          \affilnum{2}\ School of Mathematics, Shandong University, Jinan, Shandong 250100, China}

%\emails{{\tt george\_karniadakis@brown.edu} (G.~Karniadakis)}
\emails{{\tt simin\_shekarpaz@brown.edu} (S.~Shekarpaz), {\tt fanhai\_zeng@sdu.edu.cn} (F.~Zeng), {\tt george\_karniadakis@brown.edu}  
 (G.~Karniadakis)}
% \footnote and \thanks are not used in the heading section.
% Another acknowlegments/support of grants, state in Acknowledgments section
% \section*{Acknowledgments}

%%%%% Begin Abstract %%%%%%%%%%%
\begin{abstract}
We introduce a new approach for solving forward systems of differential equations using a combination of splitting methods and physics-informed neural networks (PINNs). The proposed method, splitting PINN, effectively addresses the challenge of applying PINNs to forward dynamical systems and demonstrates improved accuracy through its application to neuron models. Specifically, we apply operator splitting to decompose the original neuron model into sub-problems that are then solved using PINNs. 
Moreover, we develop an $L^1$ scheme for discretizing fractional derivatives in fractional neuron models, leading to improved accuracy and efficiency. The results of this study highlight the potential of splitting PINNs in solving both integer- and fractional-order neuron models, as well as other similar systems in computational science and engineering. 
\end{abstract}

%%%%% end %%%%%%%%%%%

%%%%% AMS/PACs/Keywords %%%%%%%%%%%

%%%% maketitle %%%%%
\maketitle

%%%% Start %%%%%%
\section{Introduction}
\label{sec1}

The human brain is a complex system that involves the interactions of billions of neurons. Mathematical models can be used to simulate the neuronal activity in the brain 
 as a system of differential equations, allowing researchers to better understand how the brain works.
Studies related to spiking neurons are performed numerically or biophysically. In numerical studies, the main goal is to solve neural equations and investigate how the dynamic behavior changes for different inputs. Biophysical approaches focus on interpreting the dynamic behavior of spiking neurons according to available experimental observations \cite{abdel2022, teka}. 

Another interesting aspect of spiking neuron models is that they can be formulated as fractional-order equations, which take into account long-term memory. The order of the derivative in these equations can affect the neuron's response \cite{Sherief2012, teka2017, weinberg2015}, making this an important area of research. Recent works in both integer- and fractional-order neuron models are discussed in Section \ref{3.3}.

In this work we introduce a new approach for solving neuron models that combines operator splitting methods with physics-informed neural networks (PINNs). Operator splitting methods have been successfully applied in various fields of physics and engineering \cite{Faou2014, Ostermann2015, shekarpaz2020, chen2019, shekarpaz2018, beck2021, li2022, haghighat2022}, while PINNs provide a powerful tool for approximating the solution of differential equations. A general introduction to the splitting method can be found in \cite{McLachlan2002ActaN2, holden2010}.

PINNs were first introduced by Raissi et al. \cite{raissi2018}. In this method, the solution of a differential equation is approximated using a neural network, and the parameters of the network are determined by solving a minimization problem that includes residual functions at collocation points, as well as initial and boundary conditions.

PINNs have been applied successfully to a broad range of ordinary and partial differential equations, including fractional equations \cite{Pang2019}, integro-differential equations, stochastic partial differential equations \cite{Zhang2020-2}, and inverse problems \cite{meng2020}. 
There have also been several extensions to the original PINN, such as Fractional PINN (FPINN) \cite{Pang2019}, physics-constrained neural networks (PCNN) \cite{liu2021, zhu2019}, variable hp-VPINN \cite{kharazmi2021}, conservative PINN (CPINN) \cite{jagtap2020-2}, Bayesian PINN \cite{yang2021}, parallel PINN \cite{raj2021}, Self-Adaptive PINN \cite{levi2020}, and Physics informed Adversarial training (PIAT) \cite{shekarpaz2022}. Innovations in activation functions, gradient optimization techniques, neural network structures, and loss function structures have driven recent advances in the field. Despite these advances, improvements are still possible, especially concerning unresolved theoretical and practical issues.

Our study makes two important contributions to the field of neural modeling. First, we propose a new method, called the splitting PINN, that employs the operator splitting technique to decompose the original spiking neuron model into sub-problems, which are then solved using PINNs. We demonstrate the effectiveness and accuracy of this method by applying it to integer- and fractional-order neuron models with oscillatory responses, for which vanilla PINN and FPINN formulations fail to predict the solutions. Second, we introduce a novel $L^1$-scheme for discretizing fractional derivatives in fractional neuron models, which leads to improved accuracy and efficiency in solving these complex models. Our results show that the combination of the splitting PINN method and the $L^1$-scheme accurately solves fractional neuron models and provides valuable insights into the underlying mechanisms of neural activity. 

This paper is organized as follows: Section 2 provides an overview of the proposed method for a given system of differential equations. Section 3 introduces various neuron models and their properties, including the Leaky Integrate-and-Fire (LIF), Izhikevich, Hodgkin-Huxley (HH), and the fractional order Hodgkin-Huxley (FO-HH) models. The efficiency and accuracy of splitting PINNs and FPINNs are demonstrated in Section 4 by applying our algorithm to different neuron models. Finally, we present a discussion of the main results in Section 5. The results and conclusion of this paper will provide valuable insights into the potential of this new approach for
solving neuron models and other similar systems in computational science and engineering.
\section{Problem setup and solution methodology}
\label{sec2}
Let us consider the general form of a nonlinear system of differential equations as follows,
\begin{equation}\label{eq2}
\frac{dx_i}{dt} = f_i(t, x), \qquad i= 1,  \cdots, n, \qquad t \in [0, T].
\end{equation}
Before introducing the splitting PINN for solving the given systems of differential equations, we briefly review
what splitting methods are, in general. 
\subsection{\textbf{Splitting method}}
For solving the given system \eqref{eq2}, let us rewrite it as follows,
$$\frac{dx}{dt} = f(t, x(t)), \qquad x(0) = x_0 \qquad x \in \mathbb{R}^n,$$
then by using splitting method, $x$ is decomposed into $(x^{*}, x^{**})$, 
where $x^{*}, x^{**}\in \mathbb{R}^d \ (d < n)$. So, we have
\begin{equation}\label{eq3}
\frac{dx^{*}}{dt} = f(t, x^{*}(t)), \qquad x^{*}(0) = x^{*}_0,
\end{equation}
\begin{equation}\label{eq4}
\frac{dx^{**}}{dt} = f(t, x^{**}(t)), \qquad x^{**}(0) = x^{**}_0,
\end{equation}
 and $x^{*}_0$ and $x^{**}_0$ are the given vectors of initial conditions. After solving the above sub-systems by the given practical algorithms, we denote the solutions of \eqref{eq3} and \eqref{eq4} as 
$$\phi_{\Delta t}^{(*)} x^{*}_0 \qquad \mathrm{and} \qquad \phi_{\Delta t}^{(**)} x^{**}_0,$$
where $\phi_{\Delta t}^{(*)}$ and $\phi_{\Delta t}^{(**)}$ are $x^{*}$-flow and $x^{**}$-flow, respectively,
and ${\Delta t}$ is the step size. Then, by combining the solutions of sub-systems, the approximation operator in $x$-flow can be written as follows 
\begin{equation}\label{eq5}
\phi_{\Delta t}^{(*)} \circ \phi_{\Delta t}^{(**)} \qquad \mbox{or} \qquad  \phi_{{\Delta t}/2}^{(**)} \circ \phi_{\Delta t}^{(*)} \circ \phi_{{\Delta t}/2}^{(**)}
\end{equation}
and $\phi_{\Delta t}^{(*)}$ and $\phi_{\Delta t}^{(**)}$ are interchangable. The first splitting is the Lie-splitting method \cite{Trotter1959OnTP}, which is first-order, and the second is the Strang splitting method \cite{Strang1968OnTC}, which is a second-order method. 

On interval $[0, T]$, we first split the original system into sub-systems (sub-problems) for each sub-interval $[t^j, t^{j+1}] 
 \ (j=0, 1, \cdots, J-1, \ t^J = T)$. Then the solution of the original system at time $t^{j+1}$ can be approximated as follows
$$x(t^{j+1}) =  \phi_{\Delta t}^{(*)} \circ \phi_{\Delta t}^{(**)} x(t^j)$$
where $x(t^j)$ is the accurate solution at $t=t^j$. 
\subsection{\textbf{Physics informed neural network}}
Consider an initial value problem as follows,
\begin{equation}\label{eq11-1}
\begin{split}
\frac{dx}{dt} &= f(t, x(t)), \qquad x \in \Omega \subseteq \mathbb{R}^n, \ t \in [0, T], \\
x(0) &= x_0,
\end{split}
\end{equation}
where $f$ is a nonlinear differential operator, and $x$ is the unknown solution with known initial condition.

By using the PINN framework, the solution of the above equation is approximated by a fully connected neural network $\mathcal{N}^L$ with $L$ layers and $N$ neurons, where the output of $l$-th hidden layer is defined as follows
\begin{equation}\label{eq7}
\mathcal{N}^l(t) = W^l \sigma (\mathcal{N}^{l-1}(t)) + b^l, \qquad 2 \leq l \leq L,
\end{equation}
$t \in \mathbb{R}$ is the input vector, $\sigma(\cdot)$ is the activation function, $\sigma (\mathcal{N}^{l-1}(t)) \in \mathbb{R}^N$  $\{W^l \in \mathbb{R}^{N\times N}, b^l \in \mathbb{R}^N\}$ are the network parameters, $\mathcal{N}^1(t) = W^1 t + b^1$ and $\mathcal{N}^L(t)$ is the output of the last layer which is used to approximate the solution. The unknown parameters can be learned by solving a minimization problem that consists of the residual error terms as follows
\begin{equation}\label{eq8}
\min_{w,b} \qquad \frac{1}{N_r} \sum_{i=1}^{N_r} |\frac{dx}{dt}(t^i) - f(t^i, x(t^i))|^2 + |x(0) - x_0|^2,
\end{equation}
where $N_r$ is the number of collocation points. The parameters are randomly initialized and optimized, and then the approximate solution is obtained.
\begin{figure}
    \centering
    \includegraphics[scale=0.40]{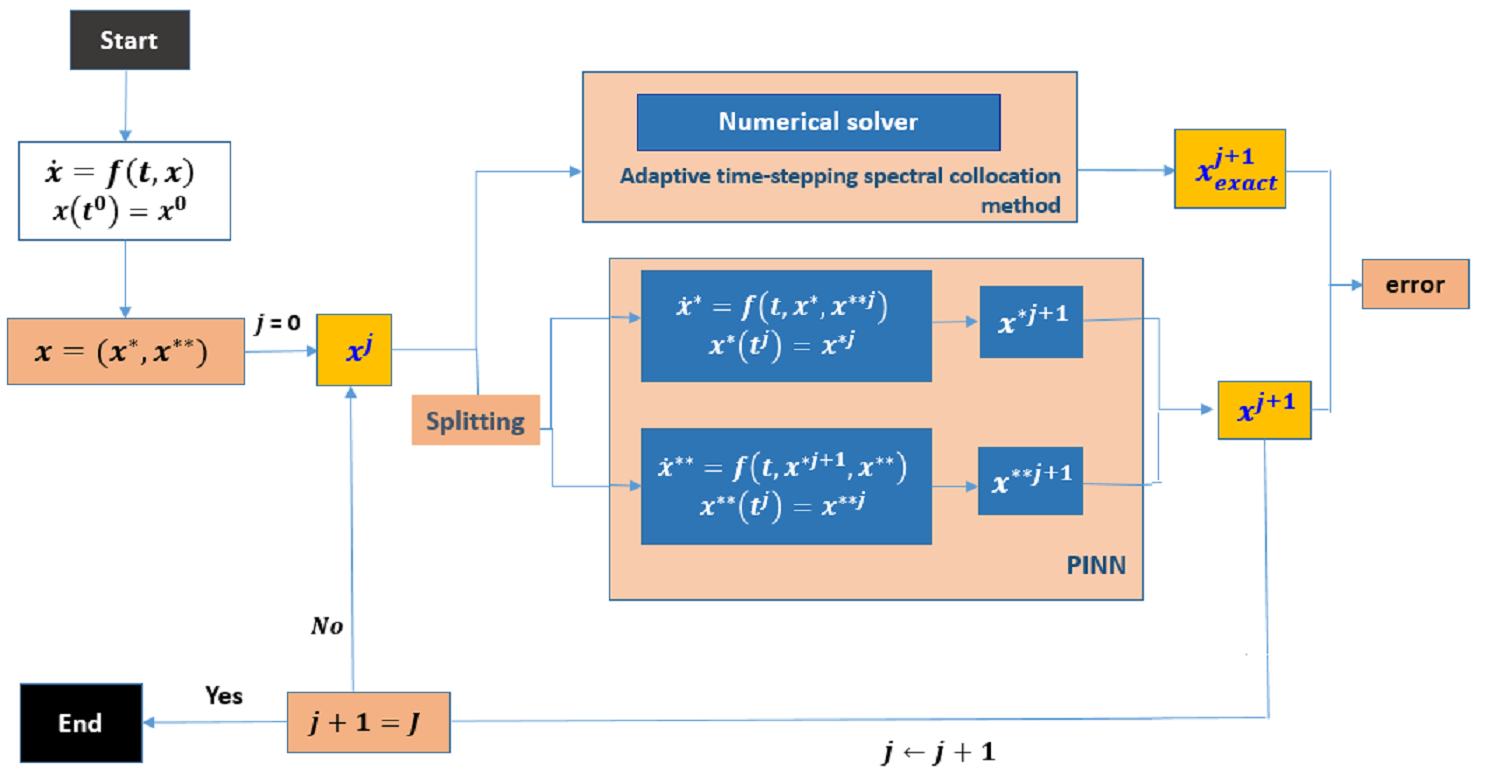}
    \caption{\textbf{Overview of the Splitting PINN :} We first split the original system into sub-systems (sub-problems) for each sub-interval $[t^j, t^{j+1}]$. For each sub-interval, $x(t^j)$ is known, the sub-systems are solved using PINN and then the solutions are combined to obtain the approximate solution $x(t^{j+1})$. To evaluate the error, we obtain the reference solution, $x_{exact}(t^{j+1})$, by using a high-order numerical solver (for more details, see Appendix B). The algorithm proceeds until arriving at a given accuracy for each sub-interval.}
    \label{1}
\end{figure}
The framework uses automatic differentiation to calculate the derivatives of the solution, which eliminates the need for manual calculation or numerical discretization \cite{AD}. This capability is available in popular deep-learning frameworks such as TensorFlow and PyTorch.
\newpage
\section{Neuron models}

Until now, various neuron models have been presented for the biological simulations of different parts of the brain. Among them, we can mention the leaky integrate-and-fire (LIF) model, the Izhikevich model, the Hodgkin-Huxley (HH) model, and FitzHugh-Nagumo (FHN) model, which model the membrane behavior \cite{abdel2022}. These models can be classified into integer-and fractional-order models. Integer-order models can capture complex phenomena in the neuron system. However, they represent only one type of firing characteristic for constant parameters of the model. On the other hand, fractional-order models can sexhibit different dynamic behavior of neurons for constant parameters \cite{teka}. This makes fractional-order models more versatile and capable of capturing a wider range of neuron behavior.

\subsection{\textbf{Neuron models: integer order}}
\subsubsection{\textbf{IF and LIF models}}
The neuron model of integrate-and-fire (IF model) is one of the best models due to the simplicity of calculations and closeness to human biological conditions. This model is a simplified version of the HH model, which is described by an equation and an assumption. Unlike the HH model, the IF model does not automatically generate an action potential. This model can be determined by 
\begin{equation}\label{eq9}
C_m \frac{dV}{dt} = I(t),
\end{equation}
with the following spike condition: if $V= V_{th}$, a spike at $t_{spike}$  is generated and the membrane potential $V(t)$ is  set to $V_{rest}$ for a refractory period $\tau_{ref}$ \cite{gerstner2002, gerstner2014}. $C_m$ is the membrane capacitance, $V_{th}$ is the voltage threshold and $V_{rest}$ is the resting membrane potential. 

A generalized type of IF model is the leaky integrate-and-fire (LIF) model, which adds a leak to the membrane potential. This model is defined by the following equation,
\begin{equation}\label{eq10}
\tau \frac{dV}{dt} = -(V-V_{rest})+RI(t),
\end{equation}
where $\tau = R C_m$ is the membrane time constant and $R$ is the membrane resistance. Because of the important properties of this model, like its computational simplicity \cite{Izhikevich2004}, accuracy in terms of the spiking behavior and spike times of neurons, and simulating speed \cite{Brette2008, Maass1996}, this model has become one of the most popular and advantageous neuron models in neuromorphic computing \cite{Aamir2018, Benjamin2014, merolla2014, Maass1996}. Also, the characteristic of membrane potential decay over time can be seen in the LIF model. 

More complex types of IF models include exponential integrate-and-fire, quadratic integrate-and-fire, and adaptive exponential integrate-and-fire \cite{Borst1999}.

\subsubsection{\textbf{Izhikevich model}}
Another neuron model for simulating the membrane behavior is the Izhikevich model. Two important features of this model are computational efficiency and biological plausibility. It reduces the more complex Hodgkin-Huxley model to a 2D system of ordinary differential equations of the form
\begin{equation}\label{eq11}
\begin{split}
&\frac{dv}{dt} = 0.04 v^2 + 5v + 140 - u + I(t), \\
&\frac{du}{dt} = a (bv - u),
\end{split}
\end{equation}
with the auxiliary condition
\begin{equation} \label{eq12}
\mbox{if} \ v \geq v_{th}, \ \mbox{then} \left \{ \begin{array}{l}
v \leftarrow c \\ 
u \leftarrow u+d
\end{array}  \right \}
\end{equation}

with $u$ and $v$ being dimensionless variables. The variable $v$ is the membrane potential of the neuron, and $u$ is the membrane recovery variable, which accounts for the activation of the $K^+$ ion current and the inactivation of the $Na^+$ ion current and provides negative feedback to the membrane potential. The dimensionless parameters $a$, $b$, $c$, and $d$ regulate the behavior of the neuron.  

The auxiliary condition \eqref{eq12} triggers a reset of the neuron when the membrane potential surpasses the threshold, simulating a spike. This model is capable of reproducing the spiking and bursting behavior of neurons in real-time, making it a widely used model in simulations of large-scale neural networks.
\\
\subsubsection{\textbf{Hodgkin-Huxley model}}
From a biophysical perspective, the nerve cell's action potential is generated by the flow of ions through the cell membrane's ion channels. Hodgkin and Huxley described the dynamics of these membrane currents through a set of coupled differential equations based on their experiments on the giant squid axon \cite{hh1952}.

\begin{figure}
    \centering
    \includegraphics[scale=0.6]{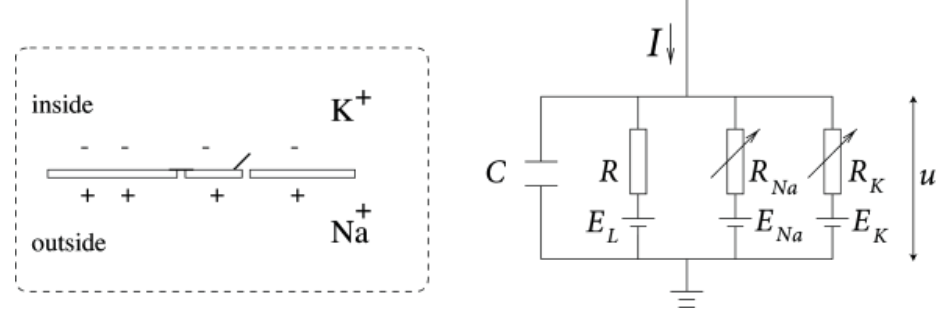}
    \caption{Schematic diagram for the Hodgkin-Huxley model  \cite{hh1952}. Left: ion channels on the membrane of the neuron. Right: simulated circuit with $I$ denoting the input current.}
    \label{3}
\end{figure}

The mechanism of the action potential can be understood with reference to Figure \ref{3}. A capacitor, resistor, and transistor were used to simulate the equivalent circuit. Changes in the action potential were observed by applying the current $I(t)$ and adjusting the capacitance and leakage resistance of the sodium and potassium channels.

  The circuit of the Hodgkin-Huxley (HH) model consists of four parallel branches: integrative branch, leaky branch, $K^+$ channel, and $Na^+$ channel. A system of four coupled differential equations was used to describe the membrane potential of a giant squid axon as follows
\begin{equation}\label{eq13-B}
\begin{split}
\frac{d V_m}{dt} =& F_1(t,V_m,n,m,h),\qquad 
\frac{d n}{dt} = F_2(t,V_m,n,m,h),\\
\frac{d m}{dt} =& F_3(t,V_m,n,m,h),\qquad
\frac{d h}{dt} = F_4(t,V_m,n,m,h),
\end{split}
\end{equation}
where 
\begin{equation}\label{eq13-C}
\begin{split}
&F_1= \frac{1}{C_m}(-g_L(V_m - E_L) - g_K n^4 (V_m - E_K) -g_{Na} m^3 h (V_m - E_{Na}) + I(t)),\\
&F_2=\alpha_n(V_m(t))(1-n(t))-\beta_n(V_m(t))n,\\
&F_3=\alpha_m(V_m(t))(1-m(t))-\beta_m(V_m(t))m,\\
&F_4=\alpha_h(V_m(t))(1-h(t))-\beta_h(V_m(t))h,
\end{split}
\end{equation}
and $g_{Na}$, $g_K$ and $g_L$ are the maximum conductances of the $Na^+$, $K^+$ and leak currents. The variables $\alpha_x$ and $\beta_x$ of the current channel conductances are dependent functions of $V_m(t)$ as follows:
\begin{equation}\label{eq14-B-2}
\begin{split}
&\alpha_n(V_m) = \frac{0.1 - 0.01(V_m - V_0)}{e^{1-0.1(V_m-V_0)} - 1}, \qquad \beta_n(V_m) = 0.125 e^{-(V_m - V_0)/80},\\
&\alpha_m(V_m) = \frac{2.5 - 0.1(V_m - V_0)}{e^{2.5-0.1(V_m-V_0)} - 1}, \qquad \beta_m(V_m) = 4.0 e^{-(V_m - V_0)/80},\\
&\alpha_h(V_m) = 0.07 e^{-(V_m-V_0)/20}, \qquad \ \ \ \ \beta_h(V_m) = \frac{1}{1+e^{3-0.1(V_m-V_0)}}.
\end{split}
\end{equation}
\subsection{\textbf{Neuron models: fractional order}}

In recent years, fractional differential equations have been developed to improve the modeling of many biological phenomena, including mechanical properties of viscoelastic tissue \cite{Magin2010}, the tissue-electrode interface
\cite{Magin2008}, pharmacokinetics of drug delivery and absorption \cite{Dokoumetzidis2009, Dokoumetzidis2010, Petras2011}, and anomalous calcium sub-diffusion in micro-domains \cite{Tan2007}.

The main characteristic of using fractional derivatives is their non-locality. This means that the next state of the system depends on the current state and all historical states before it. This advantage makes the study of fractional order systems an active area of research.
 \subsubsection{\textbf{Fractional derivative definitions}}
In this Section, we present the definitions of fractional derivatives. There are different methods for defining fractional derivatives, among which we can mention the Gr\"{u}nwald--Letnikov derivative, Riemann--Liouville derivative, and Caputo derivative \cite{fanhai2012}. The models in this paper are defined using the Caputo fractional derivative.

\begin{defn}
\label{def4}
The Caputo fractional derivative of the function $f(t)$ with order $\alpha > 0$ is defined as
\begin{equation}\label{eq200}
 \ _{C} D^{\alpha}_{a,t} f(t)  = \frac{1}{\Gamma(n-\alpha)} \int_a^t (t-s)^{n -\alpha -1} f^{(n)}(s) ds,
\end{equation} 
where $n-1 < \alpha < n$ and $n$ in a non-negative integer. 
\end{defn}

If $a=0$, then we can  use  $\frac{d^{\alpha}}{dt^{\alpha}}f(t)=\ _{C} D^{\alpha}_{0,t} f(t)$.

\subsubsection{\textbf{$L^1$ scheme to approximate the fractional derivatives}}
An efficient method for approximating the Caputo derivative of order $\alpha \ (0<\alpha<1)$ is the $L^1$ scheme, which was introduced by Oldham and Spanier \cite{Oldham1974}. 
Using this method, the Caputo derivative can be approximated by using the following formula
\begin{equation}\label{eq23}
\frac{d^{\alpha}}{dt^{\alpha}} f(t^{j}) \approx 
\delta_t^{\alpha}f^j = \sum_{k=0}^{j-1} b_{k}^{\alpha}  \left[f(t^{j-k}) - f(t^{j-1-k})\right],
\end{equation} 
where for the uniform time mesh $t^j = j \Delta t,j\ge 0$,  $\Delta t$  is the step size. Then, $b_k^{\alpha} $  is given by 
$b_k^{\alpha} = \frac{(\Delta t)^{-\alpha}}{\Gamma(2-\alpha)} [(k+1)^{1-\alpha} - k^{1-\alpha}].$
 
In \cite{lin2007}, for a smooth $f$, the error estimate of the above $L^1$ scheme  is
\begin{equation}\label{eq23-2}
 |\delta_t^{\alpha} f(t^{j}) - \frac{d^{\alpha}}{dt^{\alpha}} f(t^{j}) | \leq C \Delta t^{2-\alpha},
\end{equation} 
where $C = C(\alpha, f)$. This approximation has been used in many papers discussing fractional-order spiking neurons. The interested readers are referred to \cite{pirozzi2018, teka2018, Mondal2019, abdel2022}. Here, we will develop this method to discretize the fractional derivative, and then FPINN is used.
\subsubsection{\textbf{Fractional Hodgkin-Huxley model}}
In fractional-order neural models, the neuron's dynamics depends on the order of the derivative, which can create different types of memory-dependent dynamics. The fractional order Hodgkin-Huxley (FO-HH) model is one of the neuron models that has attracted much attention.

The HH model has two basic problems; the first is that the Dielectric losses in the membrane are neglected. The second problem is that membrane capacity is considered ideal.
To overcome the above problems, a fractional model is proposed. The idea of fractional capacity is taken from Curie's empirical law \cite{Westerlund1994}, which can be written as follows
\begin{equation}\label{eq15}
I_c(t) = C_m \frac{d^{q_1} V_m(t)}{dt^{q_1}},
\end{equation}
where $V_m(t)$ is the excitation voltage, $I_c(t)$ is the current in the capacitor, $q_1$ is the order of differentiation, and $C_m$ is the fractional capacitance \cite{weinberg20151}. The fractional order model provides a more accurate description based on long-term memory behavior. Motivated by the above discussion, we propose the following FO-HH model
\begin{equation}\label{eq16}
\begin{split}
\frac{d^{q_1} V_m}{dt^{q_1}} =& F_1(t,V_m,n,m,h),\qquad 
\frac{d^{q_2} n}{dt^{q_2}} = F_2(t,V_m,n,m,h),\\
\frac{d^{q_3} m}{dt^{q_3}} =& F_3(t,V_m,n,m,h),\qquad
\frac{d^{q_4} h}{dt^{q_4}} = F_4(t,V_m,n,m,h),
\end{split}
\end{equation}
where $q = (q_1, q_2, q_3, q_4)$ is the order of differentiation, and the other parameters are as in \eqref{eq13-B}.

By defining the fractional model, we apply the dielectric loss in the membrane, and as a result, we will see the change in the refractory period with the same value of the given current in an integer order case. The refractory period is the time when the membrane is hyperpolarized and, hence, requires a stronger stimulus to produce a smaller action potential.

The modified HH model can be very effective in biomedical applications such as heart health analysis. For example, in an ECG waveform, the PR interval represents a refractory period. PR interval estimation is very important for cardiac diagnosis \cite{Glassman1981, Chignon1993}.\\

\subsection{\textbf{Prior works in neuron models}}
\label{3.3}
In \cite{lyle2010}, the authors compared the spiking rate patterns of five single neuron models, including LIF, Izhikevich, and Hodgkin-Huxley (HH) models, under different sustained current inputs. Numerical stability and accuracy were also considered. The multi-step methods for neuronal modeling, including the HH model, were proposed in \cite{Mascagni2000}. In \cite{Khater2020}, the modified Khater (mK) method and B-spline scheme were proposed to find numerical solutions of the FitzHugh-Nagumo (FHN) equation, with a focus on finding different types of soliton wave solutions, studying their stability properties, and using them to obtain numerical solutions of the model.

In \cite{Ganguly2022}, the Hybrid Functions (HF) method was proposed as a solution for the HH model. The HF method was compared comprehensively with other algorithms, evaluating computational speed, absolute error, and integral time squared error. The finite difference scheme was used in \cite{yasin2022} to solve the stochastic FHN model, including  stability analysis and the calculation of explicit optimal a priori estimates for the existence of solutions. In \cite{Armanyos2016}, numerical solutions for the FHN and Izhikevich neuron models were obtained using a non-standard finite difference scheme and GL discretization technique. The models were compared, and their behavior was analyzed in different fractional orders.

Fractional order modeling in neural systems is a relatively new area of research. The non-local definitions of fractional calculus used in these models provide a more realistic representation of neural systems and offer a deeper understanding of their behavior. In a study reported in \cite{abdel2022}, four numerical methods were applied to two fractional-order spiking neuron models, the FO-LIF and FO-HH models. The authors used a finite memory window version of the $L^1$ approximation for comparison with well-known techniques such as the GL-based method, product integration approximation, and the Z-transform approach. In the four methods, the uniform mesh is used, and low-accurate solutions are obtained due to the singularity of the solution of fractional equations. In this paper, they have analyzed the spiking patterns, inter-spike interval adaptation, and steady-state spiking frequency for each numerical method under varying memory lengths. In a related study reported in \cite{teka}, the authors have used a $L^1$ scheme (linear interpolation based on uniform mesh) to discretize the Caputo fractional derivative. The first-order extrapolation is used to derive the linearized scheme, where the global error is $O(\Delta t^\alpha)$, and the error far from the origin is $O(\Delta t)$. In this paper, they have investigated the effects of non-Markovian power-law voltage-dependent conductances on the generation of action potentials and spiking patterns in a Hodgkin-Huxley model. They used fractional derivatives to implement the slow-adapting power-law dynamics of the potassium and sodium conductance gating variables. The results showed that, with different input currents and derivative orders, a wide range of spiking patterns can be generated, such as square wave bursting, mixed mode oscillations, and pseudo-plateau potentials. These findings suggest that power-law conductances increase the number of spiking patterns a neuron can produce. 

In \cite{Brandibur2022}, the dynamics and numerical simulations of a fractional-order coupled FHN neuronal model were discussed, and the stability properties of its equilibrium states were analyzed based on theoretical results. In \cite{tolba2019}, a non-standard finite difference scheme was used to solve the fractional Izhikevich neuron model, and a general formula for the synchronization of different Izhikevich neurons was proposed.

\section{Results}
In this Section, we use splitting PINN to solve the neuron models presented in Section 3, using the network architecture and hyperparameters specified in Table \ref{tab}. The optimization algorithm used is Adam with a learning rate of $0.0001$ and a scheduler. IN addition techniques such as adaptive activation function and feature expansion proposed in \cite{jagtap2020, lulu2022} were also used to improve the computational efficiency of the proposed method.\\
The relative $L^2$ norm of errors of the inferred  solutions are also computed as follows,
$$\mathrm{relative~L^2~error} = \frac{\sqrt {\sum_{j=1}^J (V_{exact}(t^j) - V_{app}(t^j))^2} }{\sqrt{ \sum_{j=1}^J (V_{exact}(t^j))^2}}.$$
Further details about the training procedure and hyperparameters can be found in the relevant Section.

\textbf{Data Validation:} 
In all of the examples, the reference solutions are obtained by using the adaptive time-stepping spectral collocation method, which is presented in Appendix B.

\begin{table}
\centering
	\caption{
PINN architectures and hyperparameters used for training. The first and second subcolumns under ``Depth", ``Width", and ``Activation" correspond to the first and second sub-problems, respectively.}
\label{tab}
\begin{tabular}{l|ccccc}
\hline
%\hline\noalign{\smallskip}
%\multirow{}
% \cline{1-3}
{Problems} & {Dep.} & {Wid.} & {Act.} &  {Opt.} & {Iter.}\\
\hline
\textit{LIF~model}& $5$ & $40$  & $\tanh$ & \textit{Adam} & $50000$ \\
\textit{LIF~model~with}& $7$ & $60$  & $\tanh$ & \textit{Adam} & $10000$ \\
\textit{threshold~voltage}& $$ & $$  & $$ & $$\\
\textit{Izhikevich~model}& $6, 6$ & $40, 40$  & $\tanh, \tanh$ & \textit{Adam,Adamax} & $20000$\\
\textit{HH~model}& $6, 10$ & $20, 20$  & $\tanh, \sin$ & \textit{Adam,Adamax} & $20000$\\
\textit{(step~current~function)}& $$ & $$  & $$ & $$\\
\textit{HH~model}& $6, 10$ & $20, 20$  & $\tanh, \sin$ & \textit{Adam} & $20000$\\
\textit{(constant~current)}& $$ & $$  & $$ & $$\\
\textit{FO-HH~model}& $10, 6$ & $100, 100$  & $ \tanh, \sin$ & \textit{Adam} & $70000$\\
\textit{$(q_i=0.8)$}& $$ & $$  & $$ & $$\\
\textit{FO-HH~model}& $10, 6$ & $100, 100$  & $ \tanh, \sin$ & \textit{Adam} & $50000$\\
\textit{$(q_i=0.6)$}& $$ & $$  & $$ & $$\\
\textit{FO-HH~model}& $10, 6$ & $100, 100$  & $ \tanh, \sin$ & \textit{Adam} & $20000$\\
\textit{$(q_i=0.4)$}& $$ & $$  & $$ & $$\\
\hline
\hline
\hline\noalign{\smallskip}
\end{tabular}
\end{table}

\subsection{\textbf{LIF model}}
\label{lif}
\textbf{PINN implementation:}
The approximate solutions of the LIF model can be obtained by using PINNs. For solving this model, we discretize the time domain $[0, T]$ into subdomains $[t^j, t^{j+1}]$ where $j =0, 1, \cdots, J-1$ and $t^J = T$ and then PINNs are used to solve the model in each sub-interval. More detailed visual evaluations and the numerical results of this model are provided in Appendix A.

\subsection{\textbf{Izhikevich model}}
\label{iz}
\textbf{Splitting PINN implementation:}
Consider $(u^j, v^j)$ to be the numerical solution at $t=t^j = j \Delta t$ 
 with step size, $\Delta t$, then by using the splitting method, the following sub-problems on each sub-interval $[t^j, t^{j+1}] \ (j=0, 1, \cdots, J-1, t^J = T$) are obtained
\begin{equation}\label{eq24}
\frac{du}{dt} = a (bv^j - u), \qquad u(t^j) = u^j,
\end{equation}
and
\begin{equation}\label{eq25}
\frac{dv}{dt} = 0.04 v^2 + 5v + 140 - u^{j+1} + I(t),  \qquad v(t^j) = v^{j}.
\end{equation}
This model is solved in the time interval $[0, 100 \ ms]$, with $2000$ sub-intervals and $20$ points per sub-interval. For solving sub-problems by using PINN, the fully connected neural networks are assumed to approximate the solutions, with the given architecture described in Table \ref{tab}, $v_{th} = 30mV$,
\[I(t)=\begin{cases*}
0& for $0<t<2$,\\
15 & for $t \geq 2$,
\end{cases*} \]
and the choice of parameters given in \cite{Izhikevich2003}, which can be described as follows: 
\begin{itemize}
  \item \textbf{a} $\sim 0.02 \ \frac{1}{ms}$ : Time scale of the recovery variable $u$.  Smaller values result in slower recovery.
  
  \item \textbf{b} $\sim \ 0.2 \ [dimensionless]$ : Sensitivity of the recovery variable $u$ to the sub-threshold fluctuations of the membrane potential $v$.
Greater values couple $v$ and $u$ more strongly resulting in possible sub-threshold oscillations and low-threshold spiking dynamics. The case $b<a \ (b>a)$ corresponds to saddle-node (Andronov–Hopf) bifurcation of the resting state. 

  \item \textbf{c} $\sim -50 \ mV$ : 
  The after-spike reset value of the membrane potential $v$ caused by the fast high-threshold $K^+$ conductances.
  \item \textbf{d} $\sim 2 \ mV$ :  After-spike reset of the recovery variable $u$ caused by slow high-threshold $Na^+$ and $K^+$ conductances.
\end{itemize}
\begin{figure}[!htb]
    \centering
   \begin{subfigure}[b]{0.45\textwidth}
        \includegraphics[scale=0.45]{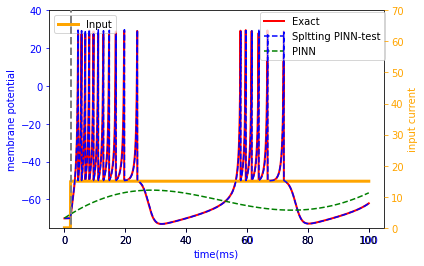}
        \caption{}
        \label{i1-06}
    \end{subfigure}
    \qquad
    \begin{subfigure}[b]{0.45\textwidth}
        \includegraphics[scale=0.45]{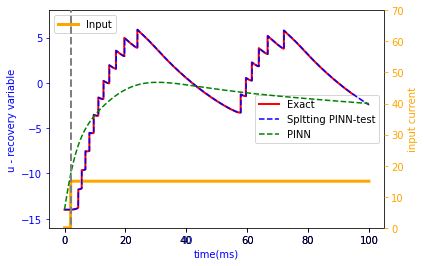}
        \caption{}
        \label{i2}
    \end{subfigure}
        \caption{Izhikevich model: Left: membrane potential versus time. Right: recovery variable versus time.}
    \label{figiz1}
        \end{figure} 

 \begin{table}
\centering
	\caption{mean relative $L^2$ norm of errors for the Izhikevich model.}
\label{tabiz1}
\begin{tabular}{lll}
\hline
%\hline\noalign{\smallskip}
%\multirow{}
% \cline{1-3}
{} & {training error} & {testing error}\\
\hline
% \textit{membrane~potential}	& $0.070808 \pm 0.003414 $   \\
\textit{membrane~potential}	& $0.115719 \pm 0.008516$  & $0.120458 \pm 0.008221$\\
\hline 
$u$  & $0.041201  \pm 0.005980$ & $0.043186 \pm 0.005070$\\
\hline
\hline\noalign{\smallskip}
\end{tabular}
\end{table}

 \begin{figure}[!htb]
    \centering
    % \begin{subfigure}[b]{0.45\textwidth}
    %     \includegraphics[scale=0.45]{err-iz.png}
    %     \caption{}
    %     \label{i2}
    % \end{subfigure}
    % \qquad
    % \begin{subfigure}[b]{0.45\textwidth}
        \includegraphics[scale=0.45]{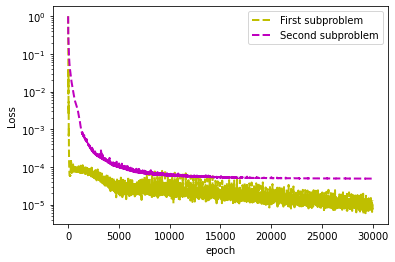}
        % \caption{}
        % \label{i1}
    % \end{subfigure}
        \caption{Loss function of the Izhikevich model.}
    \label{figiz2}
        \end{figure} 
        
\newpage
\subsection{\textbf{Hodgkin-Huxley model}}
\label{hh}
\textbf{Splitting PINN implementation:}\\
For solving the HH model by applying Splitting PINN on interval $[0, T]$, we have used the choice of parameters given in Table \ref{tabhh}.
\begin{table}
\centering
	\caption{Parameter values and their descriptions of the HH neuron model \cite{teka}.}
\label{tabhh}
\begin{tabular}{lll}
\hline
%\hline\noalign{\smallskip}
%\multirow{}
% \cline{1-3}
{Parameter} & {Value} & {Description}\\
\hline
$g_{Na}$  	& $120 \ mS/cm^2$ & \textit{Maximum~Na$^+$~current~conductance}  \\
\hline
$g_K$  & $36 \ mS/cm^2$	& \textit{Maximum~K$^+$~current~conductance} \\
\hline
$g_L$  & $0.3 \ mS/cm^2$ & \textit{Maximum~leak~current~conductance} \\
\hline
$E_{Na}$ & $50 \ mV$  & \textit{Na$^+$~current~reversal~potential} \\
\hline
$E_K$  & $-77 \ mV$ & \textit{K$^+$~current~reversal~potential} \\
\hline
$E_L$  & $-54 \ mV$ & \textit{Leak~current~reversal~potential}\\
\hline
$C_m$  & $1 \ \mu F \cdot s^{q_1-1}/cm^2$ & \textit{Membrane~Capacitance} \\
\hline
$V_0$  & $-65$ & \textit{Initial~membrane~potential} \\
\hline
$m_0$  & $0.0529$ & \textit{Initial~Na$^+$~current~activation} \\
\hline
$n_0$  & $0.3177$ & \textit{Initial~Na$^+$~current~inactivation} \\
\hline
$h_0$  & $0.5960$ & \textit{Initial~K$^+$~current~activation} \\
\hline
\hline\noalign{\smallskip}
\end{tabular}

\end{table}
Consider $(V^j, n^j, m^j, h^j)$ to be the solution at $t^j = j \Delta t$ with step size $\Delta t$. Then, the following sub-systems on each sub-interval $[t^j, t^{j+1}]$ are obtained using the Lie splitting method. The first part is,
\begin{equation}\label{eq13-B-2}
\begin{split}
&\frac{dn}{dt}=F_2(t,V_m^j,n,m,h),  \qquad n(t^j) = n^j,\qquad t\in (t^j,t^{j+1}],\\
&\frac{dm}{dt}=F_3(t,V_m^j,n,m,h),  \qquad m(t^j) = m^j,\qquad t\in (t^j,t^{j+1}],\\
& \frac{dh}{dt}=F_4(t,V_m^j,n,m,h),  \qquad h(t^j) = h^j\qquad t\in (t^j,t^{j+1}].
\end{split}
\end{equation}
where $V_m^j$ is the known solution at time $t = t^j$. The second part is  given by
\begin{equation}\label{eq14-B}
\begin{split}
&\frac{d V_m}{dt} =F_1(t, V_m, n^{j+1}, m^{j+1}, h^{j+1}) ,\qquad V_m(t^j) = V_m^j,\qquad t\in (t^j,t^{j+1}],
\end{split}
\end{equation}
and $(n^{j+1}, m^{j+1}, h^{j+1})$ is the solution obtained from the first part. The splitting approximation of the solution at $t=t^{j+1}$ is $(V_m^{j+1}, n^{j+1}, m^{j+1}, h^{j+1})$, and the solution for the whole domain is obtained by repeating this process.

Different kinds of input currents, including step current function and constant current, are considered for solving the HH model.

\begin{figure}[!htb]  
    \begin{subfigure}[b]{0.45\textwidth}    \includegraphics[scale=0.50]{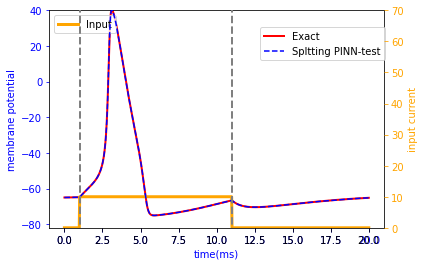}
    \caption{}
        \label{i1-1}
    \end{subfigure}
        \qquad
    \begin{subfigure}[b]{0.45\textwidth}       \includegraphics[scale=0.50]{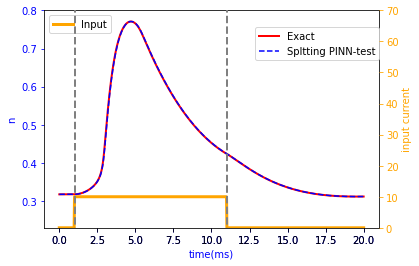}
        \caption{}
        \label{i1-2}
    \end{subfigure}
    \begin{subfigure}[b]{0.45\textwidth}
        \includegraphics[scale=0.50]{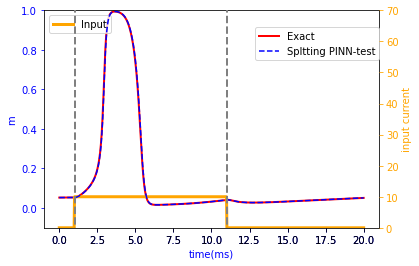}
        \caption{}
        \label{i2-0}
    \end{subfigure}
        \qquad
    \begin{subfigure}[b]{0.45\textwidth}
        \includegraphics[scale=0.50]{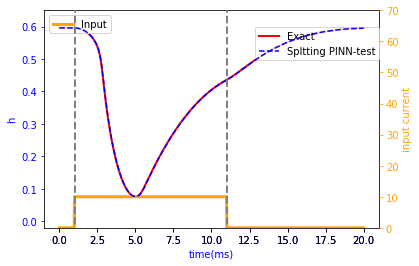}
        \caption{}
        \label{i2-01}
    \end{subfigure}
    \caption{HH model: comparison of splitting PINN results with the reference solution. (a) membrane potential; (b) activation variable of potassium channel; (c) activation variable of the sodium channel; (d) deactivation variable of the sodium channel.
    The inset plot shows the step function input current.}
    \label{fighh-step1}
    \end{figure}

\textbf{Step current function:} 
For solving the HH model with the given step current in Figure \ref{fighh-step1}(a) at the time interval $[0, 20 ms]$, we have used $800$ sub-intervals in the splitting procedure with $40$ training points for solving each sub-problem. The hyperparameters and architecture of PINNs for solving the sub-problems can be found in Table \ref{tab}. 

\textbf{Constant current:} This system is solved at time interval $[0, 100ms]$, with $3000$ sub-intervals and $30$ points at each sub-interval. With the given hyperparameters and architectures in Table \ref{tab}, PINN is used for solving each part, and the solutions are combined to obtain the approximate solutions.

The relative $L^2$ norm of errors of approximate solutions are shown in Table \ref{tabhh-fixed1}.
 \begin{table}
\centering
	\caption{mean relative $L^2$ norm of errors for HH model with current step function.}
    \label{tabhh-step1}
\begin{tabular}{lcc}
\hline
%\multirow{}
{} & {training error} & {testing error}\\
\hline
\textit{membrane~potential} & $0.001155 \pm 0.000111$  & $0.000886 \pm 0.000149$ \\

\hline
$n$ & $0.002720 \pm 5.3 \times 10^{-5}$ & $0.002712 \pm  5.3 \times 10^{-5}$ \\
\hline
$m$ & $0.014316 \pm 0.000447$ & $0.014071 \pm 0.000451 $  \\
\hline
$h$ & $0.002685 \pm 7.3 \times 10^{-5}$ & $0.002680 \pm 7.2 \times 10^{-5}$, \\
\hline 
\hline\noalign{\smallskip}
\end{tabular}
\end{table}  
\newpage
   
\begin{figure}[!htb]
    \centering
    \begin{subfigure}[b]{0.45\textwidth}      \includegraphics[scale=0.50]{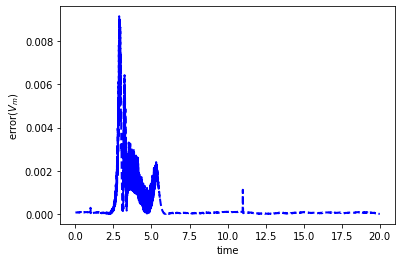}
        \caption{}
        \label{i1-3}
    \end{subfigure}
        \qquad
    \begin{subfigure}[b]{0.45\textwidth}       \includegraphics[scale=0.50]{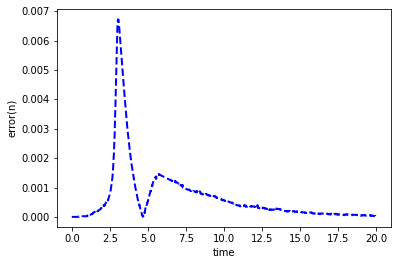}
        \caption{}
        \label{i1-4}
    \end{subfigure}
    \begin{subfigure}[b]{0.45\textwidth}
        \includegraphics[scale=0.50]{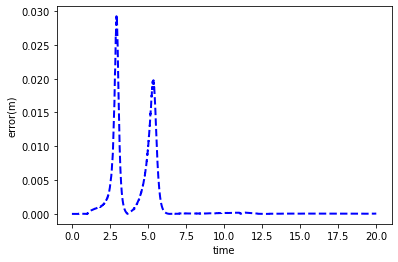}
        \caption{}
        \label{i2-0-2}
    \end{subfigure}
    \begin{subfigure}[b]{0.45\textwidth}
    \qquad        \includegraphics[scale=0.50]{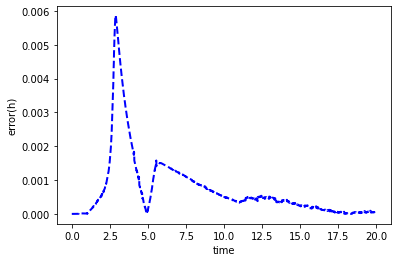}
        \caption{}
        \label{i1-5}
    \end{subfigure}
        \caption{HH model: absolute error of the normalized solutions using splitting PINN with a current step function. (a) membrane potential; (b) activation variable of potassium channel; (c) activation variable of the sodium channel; (d) deactivation variable of the sodium channel.}
    \label{fighh-step2}
        \end{figure}
        
\begin{figure}[!htb]
    \centering
   \includegraphics[scale=0.40]{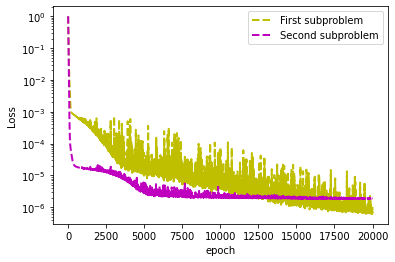}
        \caption{Loss function of the HH model with a current step function.}
    \label{fighh-step3}
        \end{figure}
$$$$        
\newpage 
 \begin{table}
\centering
	\caption{mean relative $L^2$ norm of errors for HH model for constant current.}
    \label{tabhh-fixed1}    
\begin{tabular}{lcc}
\hline
%\multirow{}
{} & {training error }& {testing error} \\
\hline
\textit{membrane~potential} & $0.009925 \pm 0.001957$ & $0.009434  \pm 0.001877 $   \\
% \hline
% \textit{membrane~potential} & $0.00162 \pm 0.001025$ & $0.003733 \pm 0.001153$   \\
\hline
$n$ & $0.005240 \pm 0.000666 $ & $0.005229 \pm 0.000661$ \\
\hline
$m$ & $0.030731 \pm 0.003931$ & $0.030364 \pm 0.003939$ \\
\hline
$h$ & $0.007077 \pm 0.000899 $ & $0.007066 \pm 0.000898 $  \\
\hline 
\hline\noalign{\smallskip}
\end{tabular}
\end{table}   
\begin{figure}[!htb]
    \centering
   \begin{subfigure}[b]{0.45\textwidth}
        \includegraphics[scale=0.50]{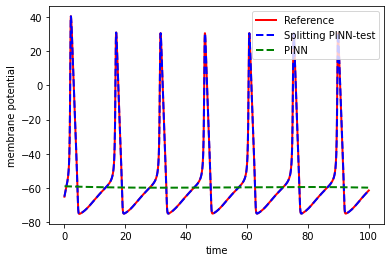}
    \caption{}
        \label{i1-6}
    \end{subfigure}
    \qquad
    ~ %add desired spacing between images, e. g. ~, \quad, \qquad, \hfill etc.
    \begin{subfigure}[b]{0.45\textwidth}
        \includegraphics[scale=0.50]{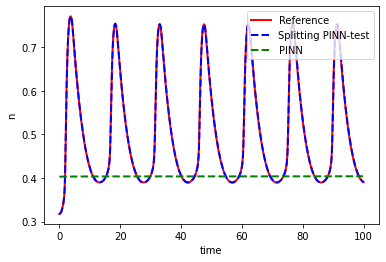}
        \caption{}
        \label{i2-0-3}
    \end{subfigure}
   \begin{subfigure}[b]{0.45\textwidth}
        \includegraphics[scale=0.50]{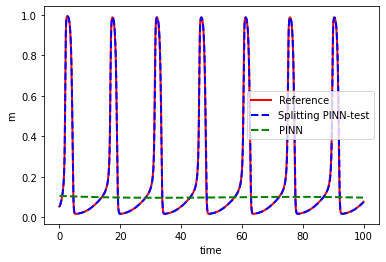}
            \caption{}
        \label{i1-7}
    \end{subfigure}
    \qquad
    ~ %add desired spacing between images, e. g. ~, \quad, \qquad, \hfill etc.
    \begin{subfigure}[b]{0.45\textwidth}
        \includegraphics[scale=0.50]{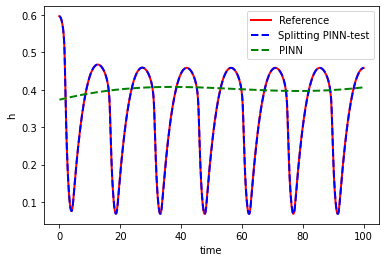}
        \caption{}
        \label{i2-0-4}
    \end{subfigure}
        \caption{HH model: comparison of splitting PINN results with the reference solution for constant current input ($I =10 \ nA/cm^2$). (a) membrane potential; (b) activation variable of potassium channel; (c) activation variable of the sodium channel; (d) deactivation variable of the sodium channel. Note that PINN fails to predict the solution.}
    \label{fighh-fixed1}
        \end{figure}

$$$$
\newpage        
 \begin{figure}[!htb]
    \centering
    \begin{subfigure}[b]{0.45\textwidth}
        \includegraphics[scale=0.45]{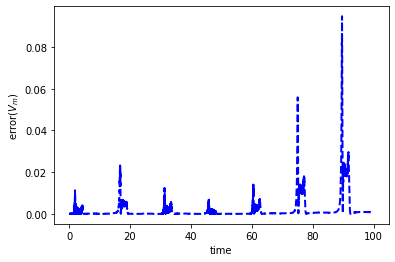}
        \caption{}
        \label{i1-8}
    \end{subfigure}
    \qquad
    \begin{subfigure}[b]{0.45\textwidth}
        \includegraphics[scale=0.45]{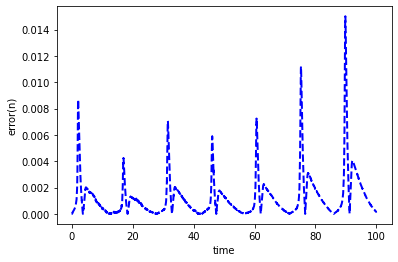}
        \caption{}
        \label{i2-0-5}
    \end{subfigure}
    \qquad
     \begin{subfigure}[b]{0.45\textwidth}
        \includegraphics[scale=0.45]{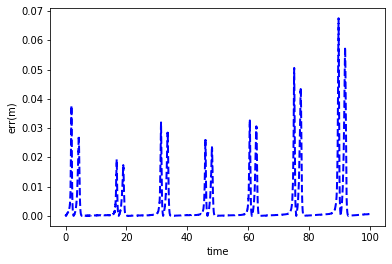}
        \caption{}
        \label{i3-3}
    \end{subfigure}
    \qquad
    \begin{subfigure}[b]{0.45\textwidth}
        \includegraphics[scale=0.45]{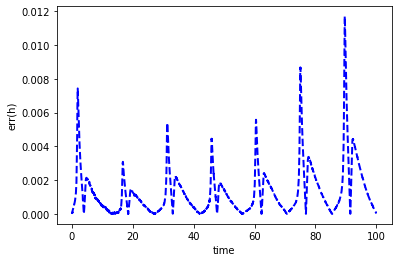}
        \caption{}
        \label{i4} 
    \end{subfigure}
        \caption{HH model: absolute error of the normalized solutions using splitting PINN with a constant current ($I =10 \ nA/cm^2$). (a) membrane potential; (b) activation variable of potassium channel; (c) activation variable of the sodium channel; (d) deactivation variable of the sodium channel.}
    \label{fighh-fixed3}
        \end{figure} 
        
\begin{figure}[!htb]
    \centering
        \includegraphics[scale=0.45]{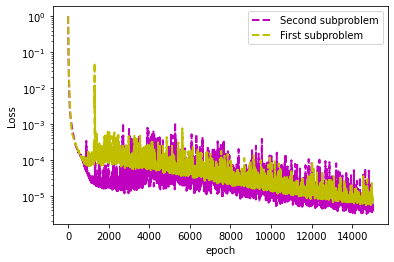}
        \caption{Loss function of the HH model for constant current.}
    \label{fighh-fixed2}
        \end{figure}

$$$$
\newpage
\subsection{\textbf{Fractional Hodgkin-Huxley model}}
\label{fhh}
\textbf{Splitting FPINN implementation:}

This problem is solved using the splitting fractional physics-informed neural network (FPINN) approach \cite{Pang2019}. In this approach, the memory of the fractional derivative at each sub-interval encompasses the entire domain. So for each sub-interval $[t^j,t^{j+1}]$ $(0 \leq j \leq N_1 - 1)$, by using the developed $L^1$ - scheme, we have the following sub-problems:

\begin{equation}\label{eq21-B-1}
\delta_t^{q_1}V_m^{N_2 j+l+1} =F_1(t^{N_2 j+l+1}, V_m^{N_2 j+l+1}, n, m, h),
\end{equation}
where 
\begin{equation}
\delta_t^{q_1}V_m^{N_2 j+l+1} = \sum_{k=0}^{N_2 j+l} b_{k}^{q_1}  \left[V_m(t^{N_2 j+l+1-k}) - V_m(t^{N_2 j+l-k})\right],
\end{equation} 
and $n$, $m$ and $h$ are the known solutions at $t = t^{N_2 j+l}$. The second subproblem is
\begin{equation}\label{eq21-B}
\begin{split}
\delta_t^{q_2}n^{N_2 j+l+1} &=F_2(t^{N_2 j+l+1}, V_m, n^{N_2 j+l+1}, m, h),\\
\delta_t^{q_3}m^{N_2 j+l+1} &=F_3(t^{N_2 j+l+1},V_m,n,m^{N_2 j+l+1},h),\\
\delta_t^{q_4}h^{N_2 j+l+1} &=F_4(t^{N_2 j+l+1},V_m,n,m,h^{N_2 j+l+1}).
\end{split}
\end{equation}
where $V_m$ is the solution of the first sub-problem at $t = t^{N_2 j+l+1}$, $l= 0, 1, \cdots, N_2 - 1$, and $N_2$ is the number of residual points for each sub-interval.

Then, the solutions of sub-systems are combined to obtain the solution of the original system.  

This model is solved with the proposed method for $q_i = 0.8, 0.6, 0.4$ at time interval $[0, 100 ms]$, with $2000$ sub-intervals and $40$ residual points in each sub-interval. Each part is solved via FPINN, with the given hyperparameters and network architectures in Table \ref{tab}.
\newpage
\begin{figure}[!htb]
    \centering
   \begin{subfigure}[b]{0.45\textwidth}
        \includegraphics[scale=0.50]{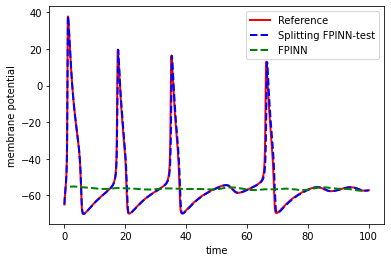}
            \caption{}
        \label{i1-9}
    \end{subfigure}
    \qquad
    \begin{subfigure}[b]{0.45\textwidth}
        \includegraphics[scale=0.50]{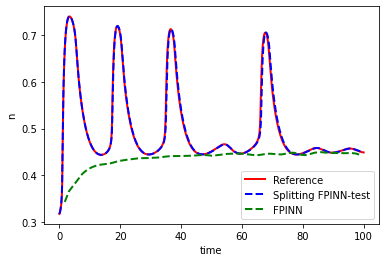}
        \caption{}
        \label{i2-0-6}
    \end{subfigure}
   \begin{subfigure}[b]{0.45\textwidth}
        \includegraphics[scale=0.50]{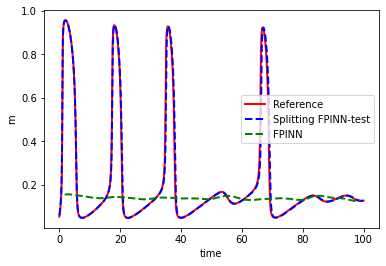}
            \caption{}
        \label{i1-10}
    \end{subfigure}
    \qquad
    \begin{subfigure}[b]{0.45\textwidth}
        \includegraphics[scale=0.50]{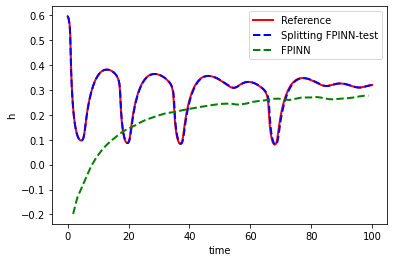}
        \caption{}
        \label{i2-0-7}
    \end{subfigure}
        \caption{
        FO-HH model: comparison of splitting FPINN results with the reference solution for constant current input ($I =20 \ nA/cm^2$) and fractional order $q_1 = q_2 = q_3 = q_4 = 0.8$. (a) membrane potential; (b) activation variable of potassium channel; (c) activation variable of the sodium channel; (d) deactivation variable of the sodium channel. Note that FPINN fails to predict the solution.
        }
\label{figfhh-0.8-1}
        \end{figure}

\begin{figure}[!htb]
    \centering
   % \begin{subfigure}[b]{0.45\textwidth}
        \includegraphics[scale=0.45]{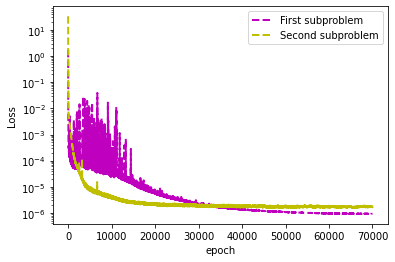}
        \label{i1-11}
  \caption{Loss function of the FO-HH model for fractional order  $q_1 = q_2 = q_3 = q_4 = 0.8$.}
\label{figfhh-0.8-2}
        \end{figure} 

\begin{figure}[!htb]
    \centering
   \begin{subfigure}[b]{0.45\textwidth}
        \includegraphics[scale=0.50]{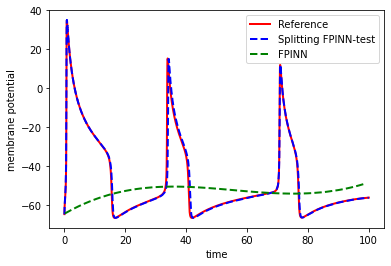}
            \caption{}
        \label{i1-12}
    \end{subfigure}
    \qquad
    \begin{subfigure}[b]{0.45\textwidth}
        \includegraphics[scale=0.50]{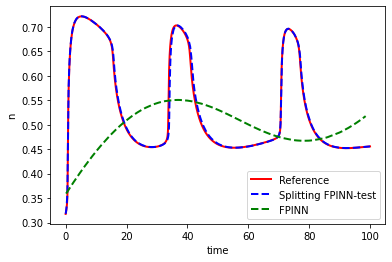}
        \caption{}
        \label{i2-0-8}
    \end{subfigure}
   \begin{subfigure}[b]{0.45\textwidth}
        \includegraphics[scale=0.50]{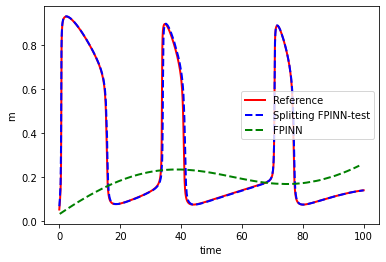}
            \caption{}
        \label{i1-13}
    \end{subfigure}
    \qquad
    \begin{subfigure}[b]{0.45\textwidth}
        \includegraphics[scale=0.50]{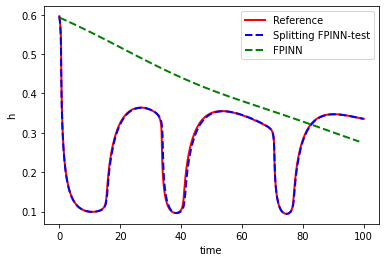}
        \caption{}
        \label{i2-0-9}
    \end{subfigure}
        \caption{FO-HH model: comparison of splitting FPINN results with the reference solution for constant current input ($I =20 \ nA/cm^2$) and fractional order $q_1 = q_2 = q_3 = q_4 = 0.6$. (a) membrane potential; (b) activation variable of potassium channel; (c) activation variable of the sodium channel; (d) deactivation variable of the sodium channel. Note that FPINN fails to predict the solution.}
\label{figfhh-0.6-1}
        \end{figure}

\begin{figure}[!htb]
    \centering
   % \begin{subfigure}[b]{0.45\textwidth}
        \includegraphics[scale=0.45]{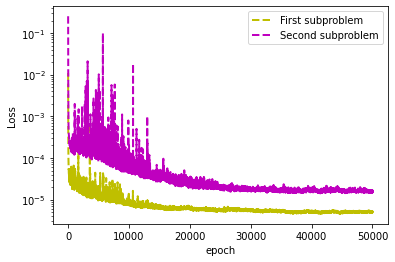}
        \label{i1-14}
  \caption{Loss function of the FO-HH model for fractional order  $q_1 = q_2 = q_3 = q_4 = 0.6$.}
\label{figfhh-0.6-2}
        \end{figure}         
  
\begin{figure}[!htb]
    \centering
   \begin{subfigure}[b]{0.45\textwidth}
        \includegraphics[scale=0.50]{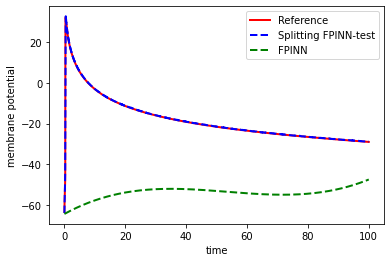}
            \caption{}
        \label{i1-15}
    \end{subfigure}
    \qquad
    \begin{subfigure}[b]{0.45\textwidth}
        \includegraphics[scale=0.50]{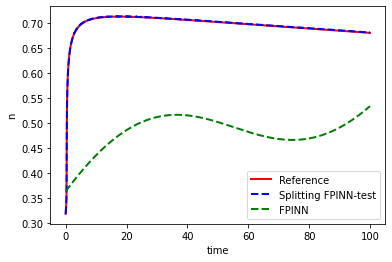}
        \caption{}
        \label{i2-0-10}
    \end{subfigure}
   \begin{subfigure}[b]{0.45\textwidth}
        \includegraphics[scale=0.50]{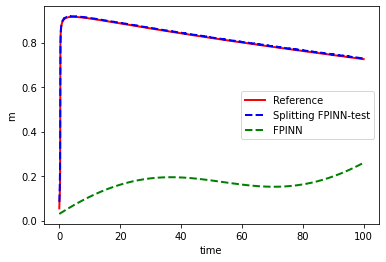}
            \caption{}
        \label{i1-16}
    \end{subfigure}
    \qquad
    \begin{subfigure}[b]{0.45\textwidth}
        \includegraphics[scale=0.50]{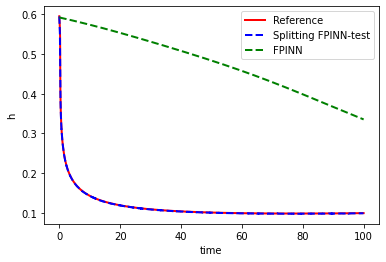}
        \caption{}
        \label{i2-0-11}
    \end{subfigure}
        \caption{FO-HH model: comparison of splitting FPINN results with the reference solution for constant current input ($I =20 \ nA/cm^2$) and fractional order $q_1 = q_2 = q_3 = q_4 = 0.4$. (a) membrane potential; (b) activation variable of potassium channel; (c) activation variable of the sodium channel; (d) deactivation variable of the sodium channel. Note that FPINN fails to predict the solution.}
        \label{figfhh-0.4-1}
        \end{figure}
 % \begin{figure}[!htb]
 %    \centering
 %    \begin{subfigure}[b]{0.45\textwidth}
 %        \includegraphics[scale=0.45]{errv-hh0.4-.png}
 %        \caption{}
 %        \label{i1}
 %    \end{subfigure}
 %    \qquad
 %    \begin{subfigure}[b]{0.45\textwidth}
 %        \includegraphics[scale=0.45]{errn-fhh0.4.png}
 %        \caption{}
 %        \label{i2}
 %    \end{subfigure}
 %    \qquad
 %     \begin{subfigure}[b]{0.45\textwidth}
 %        \includegraphics[scale=0.45]{errm-fhh0.4.png}
 %        \caption{}
 %        \label{i3}
 %    \end{subfigure}
 %    \qquad
 %    \begin{subfigure}[b]{0.45\textwidth}
 %        \includegraphics[scale=0.45]{errh-fhh0.4.png}
 %        \caption{}
 %        \label{i4} 
 %    \end{subfigure}
 %        \caption{Absolute error of normalized solutions for the FHH model with $q_i = 0.4$}
 %    \label{fighh-fixed3}
 %        \end{figure} 
        
\begin{figure}[!htb]
    \centering
   % \begin{subfigure}[b]{0.45\textwidth}
        \includegraphics[scale=0.45]{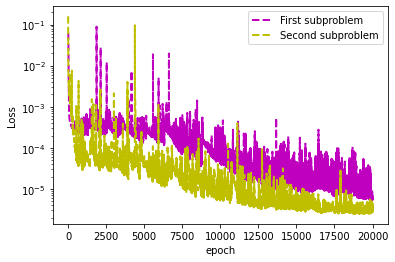}
        \label{i1-17}
  \caption{Loss function of the FO-HH model for fractional order $q_1 = q_2 = q_3 = q_4 = 0.4$.}
\label{figfhh-0.4-2}
        \end{figure} 
        
 \begin{table}
\centering
	\caption{mean relative $L^2$ norm of errors for FO-HH model for different orders of fractional derivatives. The first and
second sub-rows correspond
to the training and testing errors, respectively.}
\label{tabfhh-1}
\begin{tabular}{lccc}
\hline
%\hline\noalign{\smallskip}
%\multirow{}
% \cline{1-3}
{} & {$q_i = 0.8$} & {$q_i = 0.6$} & {$q_i = 0.4$}\\
\hline
\textit{membrane}	& $0.084252 \pm 0.011138$ & $0.098086 \pm 0.027659$ & $0.043236 \pm 0.001514$  \\
\textit{potential}	& $0.084460 \pm 0.012576$ & $0.096949 \pm 0.028145$ & $0.043241 \pm 0.003771$  \\
\hline
$n$  & $0.024353 \pm 0.005469$ & $0.023854 \pm 0.008417$ & $0.001909 \pm 8.8 \times 10^{-5}$\\
\textit{}	& $0.024360 \pm 0.005527$ & $0.023776 \pm 0.008435$ & $0.001933 \pm 0.000177$  \\
\hline
$m$  & $0.129660 \pm 0.035021$ & $0.093863 \pm 0.031778$ & $0.006009 \pm 9.4 \times 10^{-5} $\\
\textit{}	& $0.129083 \pm 0.035903$ & $0.093619 \pm  0.032427$ & $0.006063 \pm 0.000214$  \\
\hline
$h$ & $0.036292 \pm 0.009040$ &  $0.040947 \pm 0.013980$ & $0.007536 \pm 0.000374$\\
\textit{}	& $0.036336 \pm 0.009121$ & $0.040921 \pm 0.013963$ & $0.007592 \pm 0.000649$  \\
\hline
\hline\noalign{\smallskip}
\end{tabular}

\end{table}

\begin{figure}[!htb]
    \centering
   \begin{subfigure}[b]{0.45\textwidth}
        \includegraphics[scale=0.50]{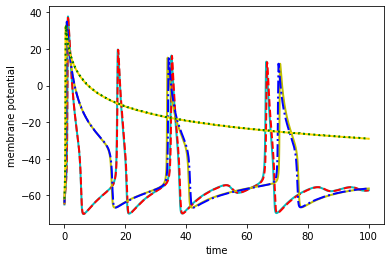}
            \caption{}
        \label{i1-18}
    \end{subfigure}
    \qquad
    \begin{subfigure}[b]{0.45\textwidth}
        \includegraphics[scale=0.50]{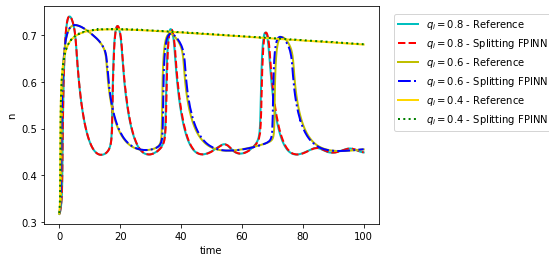}
        \caption{}
        \label{i2-0-12}
    \end{subfigure}
   \begin{subfigure}[b]{0.45\textwidth}
        \includegraphics[scale=0.50]{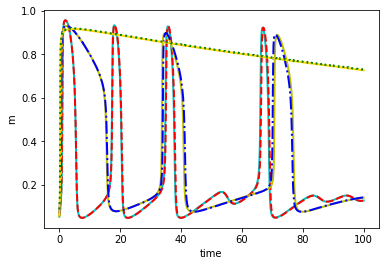}
            \caption{}
        \label{i1-19}
    \end{subfigure}
    \qquad
    \begin{subfigure}[b]{0.45\textwidth}
        \includegraphics[scale=0.50]{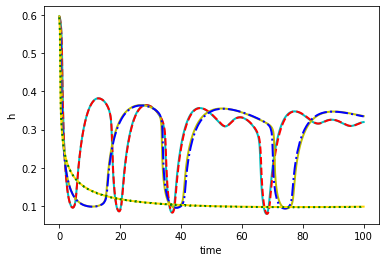}
        \caption{}
        \label{i2-0-13}
    \end{subfigure}
        \caption{Plots of the membrane potential and the corresponding activation and deactivation variables of the FO-HH model with constant current $I=20 \ nA/cm^2$ for different orders of fractional derivatives, obtained using the splitting FPINN.}
    \label{figfhh}
        \end{figure}
        
\newpage
\section{Disccusion}
We have presented a deep learning approach for solving nonlinear systems of differential equations. The performance and accuracy of splitting Physics-Informed Neural Networks (PINNs) are studied in the context of solving neuron models. In Section \ref{iz}, the proposed method was used to solve the Izhikevich model. We showed the reference and learned solutions in Figure \ref{figiz1} for equi-distant training and random test datasets. It can be seen that the nonlinear solutions learned by the nonlinear splitting PINN match well with the reference solutions. The relative $L^2$ norm of errors for the solutions are presented in Table \ref{tabiz1}, which shows the method's capability for this model.
The convergence of the loss functions is shown in Figure \ref{figiz2}. It can be seen that this function rapidly tends to values lower than $10^{-5}$ and $10^{-4}$ for the first and second sub-problems. 

After solving the Izhikevich model, we focus on using the Splitting PINN method for the Hodgkin-Huxley model. In this model, applying input current and changes in voltage resulting from the opening and closing of ion channels lead to the generation of spikes. The amount of input current needed to generate a spike is at least $2.7 nA/cm^2$, and with increasing input current, the number of spikes increases. In Section \ref{hh}, the voltage is obtained for two states of input current:  step function current and constant input current. In the first case, when the current is applied over a given time interval, the voltage increases and produces an action potential (positive peak). After the spike, the potassium channel opens, and the sodium channel closes, causing the voltage to decrease and the neuron to enter a refractory period where the potential is below the resting potential of $V_{rest} = -65 mV$. The voltage then slowly returns to $-65 mV$. The splitting PINN method was used to calculate approximate solutions and compare them with reference solutions. The results are shown in Figure \ref{fighh-step1}, and the absolute errors are given in Figure \ref{fighh-step2}, which demonstrate the accuracy of the method. The relative $L^2$ norm of errors for the current step function is presented in Table \ref{tabhh-step1}. The convergence of the loss functions is also shown in Figure \ref{fighh-step3} and tends to smaller values of $10^{-6}$ and $10^{-5}$ for the first and second sub-problems, respectively.

In the case of a constant input current, the approximate and reference solutions are shown in Figure \ref{fighh-fixed1}. The proposed method's effectiveness is demonstrated by comparing the approximate results with the reference and numerical solutions obtained using PINN. The absolute errors of the solutions are shown in Figure \ref{fighh-fixed3}, and the plots of the loss functions are displayed in Figure \ref{fighh-fixed2}, which converge to values of $10^{-6}$ for sub-problems. The relative $L^2$  norm of errors of the solutions is presented in Table \ref{tabhh-fixed1}. These results demonstrate the capability and accuracy of the proposed algorithm in solving this model.
    
The fractional-order Hodgkin-Huxley model was also studied in Section \ref{fhh}. The system was modeled using Caputo's fractional derivative and simulated using the proposed method. The effect of memory, associated with fractional order, on firing activity was also examined. To address this, the $L^1$ scheme was improved for solving the fractional-order model by proposing using full-domain memory, rather than sub-domain memory, to calculate the fractional derivative. The numerical results obtained using the improved method for $q_i = 0.8, 0.6, 0.4$ are shown in Figures \ref{figfhh-0.8-1}, \ref{figfhh-0.6-1} and \ref{figfhh-0.4-1}, where the approximate results are compared with the reference solutions. The plots of the loss functions are also shown in Figures \ref{figfhh-0.8-2}, \ref{figfhh-0.6-2}, and \ref{figfhh-0.4-2}, demonstrating the efficiency and convergence of the current method. The relative $L^2$ norm of errors of the solutions is presented in Table \ref{tabfhh-1}, and as shown in the Section, the results are accurate and have converged to the reference solutions.

The voltage responses for the different orders of fractional derivatives under constant input current $I=20$ are also displayed in Figure \ref{figfhh}. The Figure demonstrates that various spike patterns can be produced. As $q_i$ approaches 1, the firing frequency decreases and results in an increase in the number of spikes in the same time period. Furthermore, the first spike occurs at a later time. Following the cessation of the injected current, the memory-dependent spiking activity can also be observed by applying a step function current.
Additionally, the regularity of solutions is observed to depend on the order of fractional derivatives in Figure \ref{figfhh}. As $q_i$ decreases, the regularity increases. Table \ref{tabfhh-1} shows that, unregular solutions are more sensitive to parameter initialization, leading to a greater standard deviation. To address this issue, the network architecture can be improved. 

The results in this section demonstrate high accuracy, and the solutions have converged to the reference values.

\section{\textbf{Conclusion}}
In this study, we introduced a novel method for solving neuron models represented as systems of differential equations. The Splitting PINN algorithm was demonstrated to be effective and accurate through comparisons with reference solutions and the mean $L^2$ relative norm of errors. In addition, the results of the fractional-order Hodgkin-Huxley model highlighted the effect of memory on firing activity and voltage responses, showing that as $q_i$ approaches 1, the firing frequency decreases while the number of spikes in the same time period increases.

This research provides valuable insights into the behavior of neuron membranes and the various spike patterns that can be generated. The performance of the proposed method demonstrates its superiority over the vanilla PINN algorithm in solving complex neuron models represented as systems of differential equations. This study contributes to the development of tools for investigating the behavior of neurons and the underlying mechanisms of neural activity.

In conclusion, the proposed Splitting PINN algorithm is a promising method for solving systems of differential equations and provides valuable insights into the behavior of neurons and their underlying mechanisms.

%%%% Acknowledgments %%%%%%%%
\section*{Acknowledgments}
We would like to thank Dr. Khemraj Shukla and Dr. Ehsan Kharazmi for helpful discussions. This work was supported by AFOSR MURI funding (FA9550-20-1-0358) and Graphs and Spikes for Earth and
Embedded Systems (SEA-CROGS) project. Fanhai Zeng is supported by the National Natural Science Foundation of China (12171283),
the National Key R\&D Program of China (2021YFA1000202, 2021YFA1000200), 
the Science Foundation Program for Distinguished Young Scholars of Shandong (Overseas) (2022HWYQ-045).

\appendix

\newpage
% \appendix
% \section{PINN implementation for LIF model}
% \label{appen1}
%\textbf{Appendix A: PINN implementation for LIF model}\\

\section{PINN implementation for LIF model}

% \textbf{Validation:} The reference solution of the LIF model is obtained by using adaptive time-stepping spectral collocation method \cite{Zhang2020}.\\

\begin{figure}[!htb]
    \centering
   \begin{subfigure}[b]{0.38\textwidth}
        \includegraphics[scale=0.44]{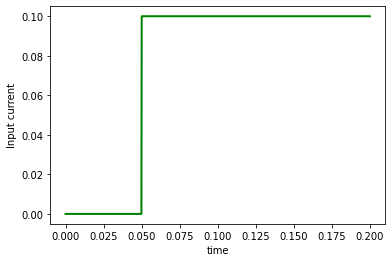}
        \caption{$I(t)$}
        \label{i1-00}
    \end{subfigure}
    \qquad
        \qquad
   \begin{subfigure}[b]{0.38\textwidth}
        \includegraphics[scale=0.44]{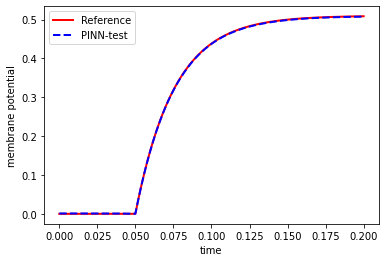}
        \caption{$V(t)$}
        \label{i1-01}
    \end{subfigure}
        \caption{LIF model:  Left: Input current versus time. Right: membrane potential versus time (comparison of PINN results with the reference solutions).}
    \label{Fig04}
        \end{figure}

In this Section, PINN has been used to solve the LIF model with the given parameters in Table \ref{tab2.2}. For different input currents in time interval $[0, 0.5s]$, the network is trained with $1000$ residual points. Initially, assuming that there is no input current, then if the current $I(t) =0.1 A$ is injected into the neuron, the membrane potential tends to $R \cdot I(t)$, which can be seen in Figure \ref{Fig04}.

\begin{table}
\centering
	\caption{Parameter values and their descriptions of LIF model \cite{eshraghian2021-2, eshraghian2021}.}
\label{tab2.2}
\begin{tabular}{lll}
\hline
%\hline\noalign{\smallskip}
%\multirow{}
% \cline{1-3}
{Parameter} & {Value} & {Description}\\
\hline
$C_m$  	& $5\times 10^{-3} F$ & \textit{Membrane~Capacity}  \\
\hline
$R$  & $5.1 \ \Omega$	& \textit{Membrane~Resistance} \\
\hline
$\tau$  & $RC_m \ s$ & \textit{Membrane~time~constant} \\
\hline
$V_{rest}$ & $ 0 \ V$  & \textit{Resting~membrane~potential} \\
\hline
$V_{th}$  & $ 1 \ V$ & \textit{Threshold~membrane~potential} \\
\hline
\hline\noalign{\smallskip}
\end{tabular}
\end{table}

\begin{figure}[!htb]
    \centering
   \begin{subfigure}[b]{0.28\textwidth}
   \includegraphics[scale=0.33]{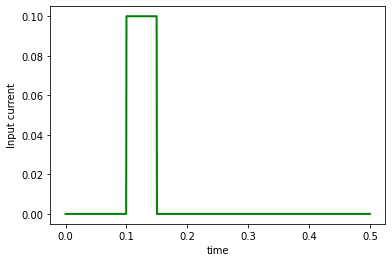}
        \caption{$I_1(t)$}
        \label{i102}
    \end{subfigure}
    \qquad
    \begin{subfigure}[b]{0.24\textwidth}
        \includegraphics[scale=0.33]{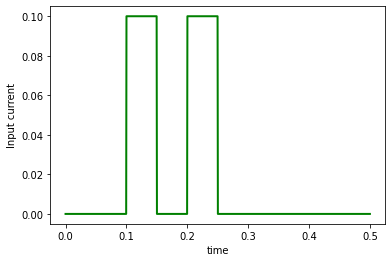}
        \caption{$I_2(t)$}
        \label{i2-015}
    \end{subfigure}
    \qquad
    ~ %add desired spacing between images, e. g. ~, \quad, \qquad, \hfill etc.
    \begin{subfigure}[b]{0.24\textwidth}
        \includegraphics[scale=0.33]{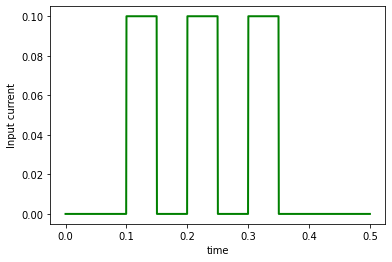}
        \caption{$I_3(t)$}
        \label{i3-1}
    \end{subfigure}    

  \begin{subfigure}[b]{0.28\textwidth}
        \includegraphics[scale=0.33]{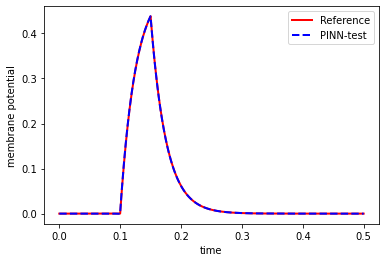}
        \caption{$V_1(t)$}
        \label{u1-4}
    \end{subfigure}
    \qquad
    \begin{subfigure}[b]{0.28\textwidth}
        \includegraphics[scale=0.33]{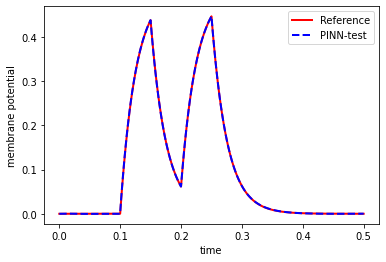}
        \caption{$V_2(t)$}
        \label{u2-8}
    \end{subfigure}
    \qquad
    \begin{subfigure}[b]{0.28\textwidth}
        \includegraphics[scale=0.33]{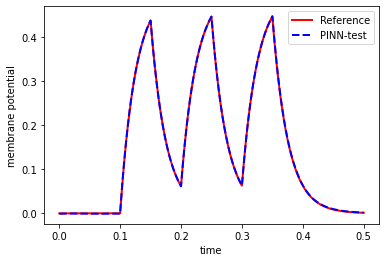}
        \caption{$V_3(t)$}
        \label{u2-1}
    \end{subfigure} 
    \begin{subfigure}[b]{0.28\textwidth}
        \includegraphics[scale=0.33]{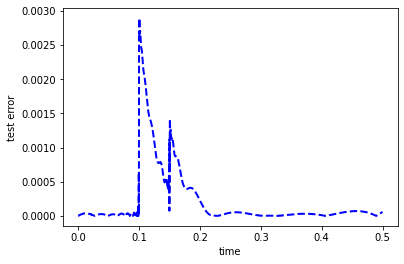}
        \caption{absolute error of $V_1(t)$}
        \label{u1-1}
    \end{subfigure}
    \qquad
    \begin{subfigure}[b]{0.28\textwidth}
        \includegraphics[scale=0.33]{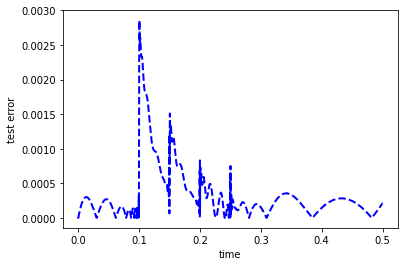}
        \caption{absolute error of $V_2(t)$}
        \label{u2-2}
    \end{subfigure}
    \qquad
    \begin{subfigure}[b]{0.28\textwidth}
        \includegraphics[scale=0.33]{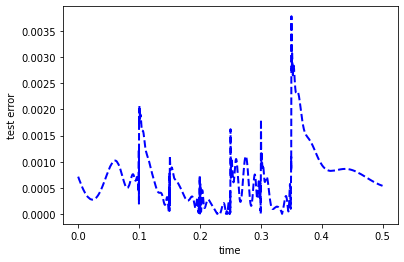}
        \caption{absolute error of $V_3(t)$}
        \label{u2-3}
    \end{subfigure}      
        \caption{
        LIF model; (a)-(c): input step function currents, (d)-(f): comparison of the PINN results with the reference solutions, and (g)-(i): the error results of PINN.}
    \label{Fig1-1}
        \end{figure}
     
\begin{figure}[!htb]
    \centering
   \begin{subfigure}[b]{0.28\textwidth}
        \includegraphics[scale=0.33]{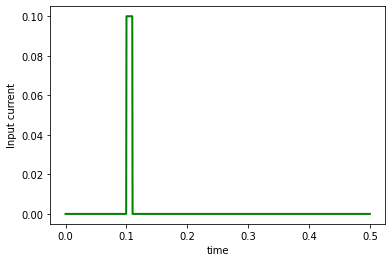}
        \caption{$I_1(t)$}
        \label{ti5}
    \end{subfigure}
    \qquad
    \begin{subfigure}[b]{0.28\textwidth}
        \includegraphics[scale=0.33]{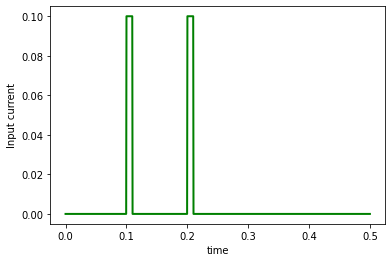}
        \caption{$I_2(t)$}
        \label{i2-1}
    \end{subfigure}
    \qquad
    \begin{subfigure}[b]{0.28\textwidth}
        \includegraphics[scale=0.33]{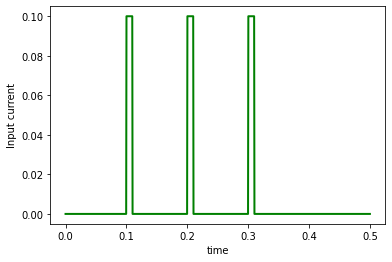}
        \caption{$I_3(t)$}
        \label{i3-2}
    \end{subfigure}   

  \begin{subfigure}[b]{0.28\textwidth}
        \includegraphics[scale=0.33]{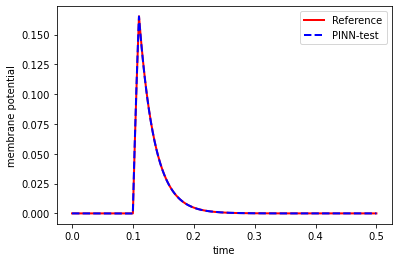}
        \caption{$V_1(t)$}
        \label{u1-2}
    \end{subfigure}
    \qquad
    \begin{subfigure}[b]{0.28\textwidth}
        \includegraphics[scale=0.33]{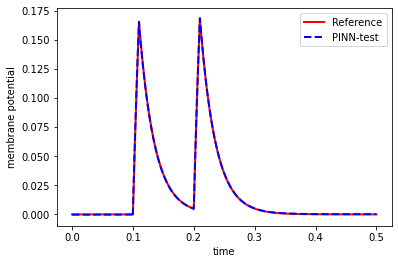}
        \caption{$V_2(t)$}
        \label{u2-4}
    \end{subfigure}
    \qquad
    \begin{subfigure}[b]{0.28\textwidth}
        \includegraphics[scale=0.33]{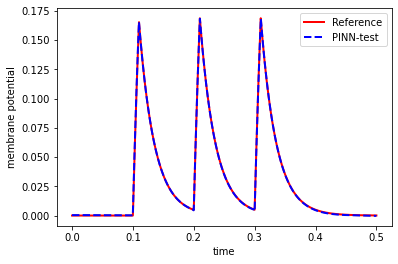}
        \caption{$V_3(t)$}
        \label{u2-5}
    \end{subfigure} 

    \begin{subfigure}[b]{0.28\textwidth}
        \includegraphics[scale=0.33]{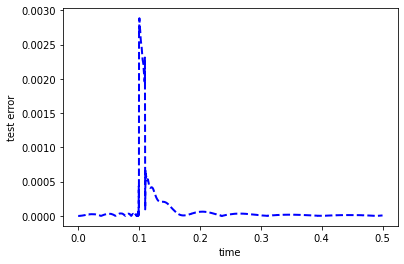}
        \caption{absolute error of $V_1(t)$}
        \label{u1-3}
    \end{subfigure}
    \qquad
    \begin{subfigure}[b]{0.28\textwidth}
        \includegraphics[scale=0.33]{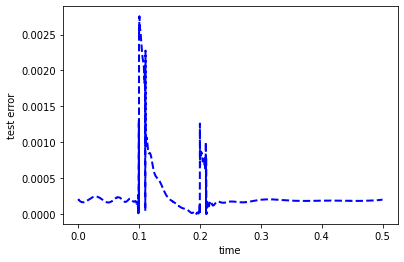}
        \caption{absolute error of $V_2(t)$}
        \label{u2-6}
    \end{subfigure}
    \qquad
    \begin{subfigure}[b]{0.28\textwidth}
        \includegraphics[scale=0.33]{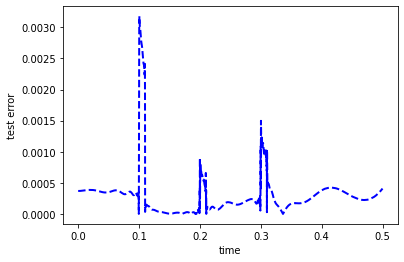}
        \caption{absolute error of $V_3(t)$}
        \label{u2-7}
    \end{subfigure}      
        \caption{
        LIF model; (a)-(c): input spikes, (d)-(f): comparison of the PINN results with the reference solutions, and (g)-(i): the error results of PINN.}
    \label{Fig2-1}
        \end{figure}
        
\begin{figure}[!htb]
    \centering
   \begin{subfigure}[b]{0.39\textwidth}
        \includegraphics[scale=0.39]{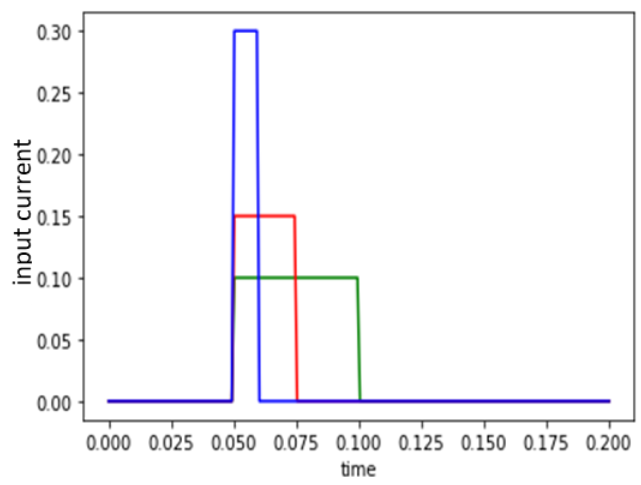}
        \caption{different input currents}
        \label{i1-03}
    \end{subfigure}
    \qquad
    ~ %add desired spacing between images, e. g. ~, \quad, \qquad, \hfill etc.
    \begin{subfigure}[b]{0.39\textwidth}
        \includegraphics[scale=0.39]{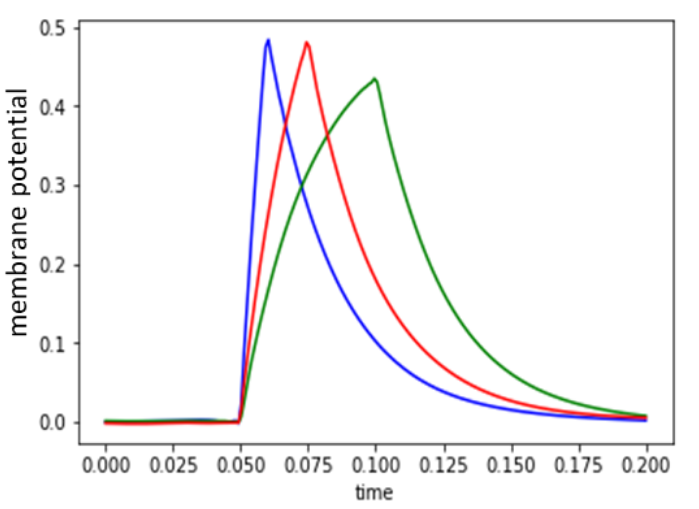}
        \caption{membrane potentials}
        \label{i2-2}
    \end{subfigure}
        \caption{LIF model: (a) Input currents with different lengths and heights, (b) corresponding membrane potentials by using PINN}
    \label{Fig3}
        \end{figure}

The membrane potentials for different forms of step current function have been displayed in Figures \ref{Fig1-1}--\ref{Fig3}. By applying an input current in a time step, the membrane potential increases, and then after cutting off the input current, the potential decreases with the time constant $\tau$.

If the width of the injected current increases, the membrane potential increases at a slower rate. As the input current pulse amplitude decreases, the membrane potential jumps directly up in a short time. The membrane potential decreases over time in the absence of an input current.

\newpage
So far, we have seen the neuron's behavior against different input currents (spikes). We need to apply a voltage threshold condition to the model to see the process of generating output spikes by a neuron. According to this condition, if the membrane potential exceeds this threshold, the neuron produces a spike, and after firing the spike, the membrane potential returns to its initial state (resting state).

If the injected current increases, the membrane potential approaches the threshold faster, and the firing rate increases. Similar behavior in firing frequency can be produced by lowering the threshold. These experiments hold for the constant current application, but neurons may be more likely to receive spikes in other cases.

\begin{figure}[!htb]
    \centering
   \begin{subfigure}[b]{0.25\textwidth}
        \includegraphics[scale=0.31]{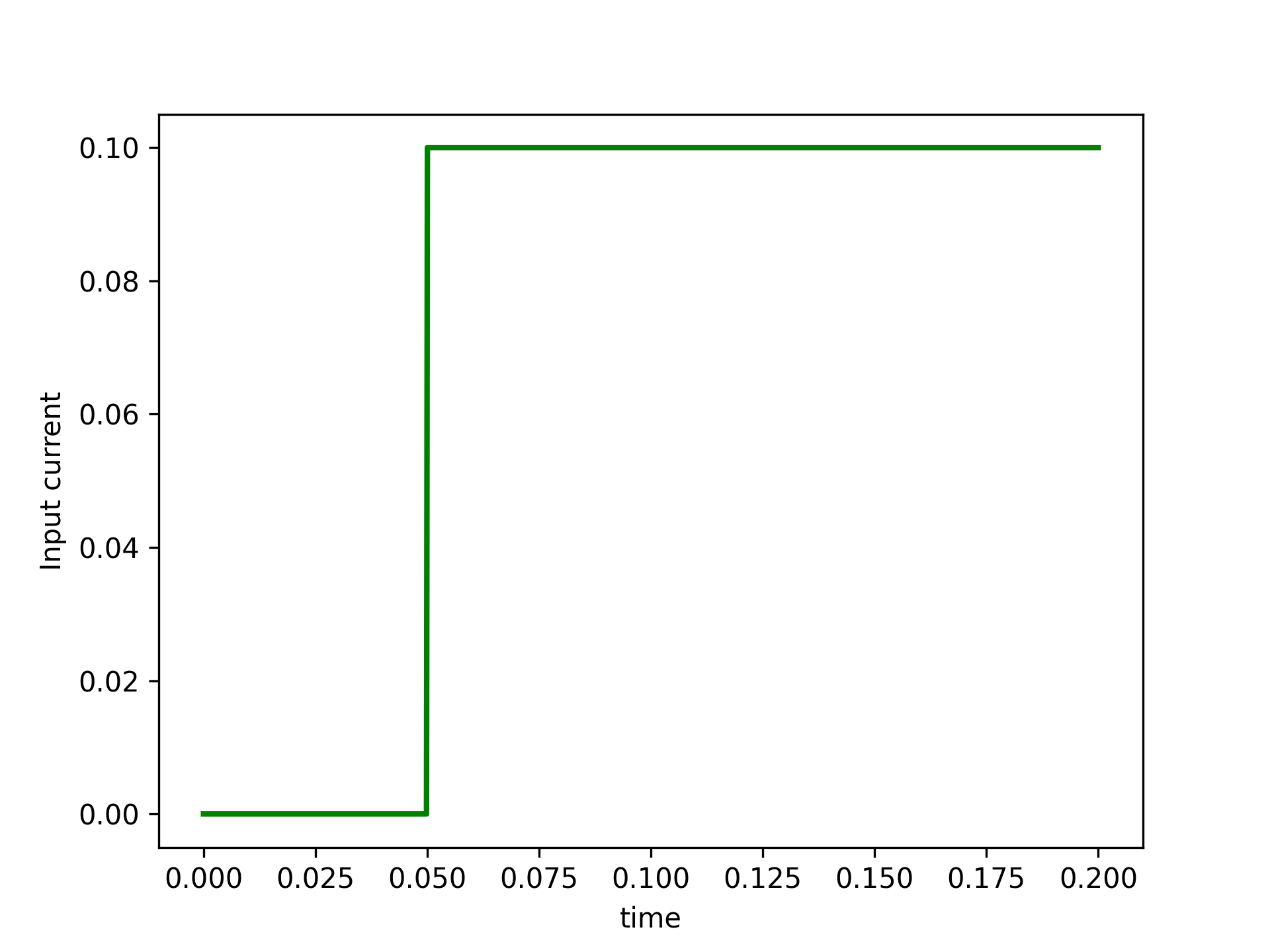}
        \caption{input current $I_{1}(t)$}
        \label{i1-0}
    \end{subfigure}
    \qquad
    \qquad
   \begin{subfigure}[b]{0.25\textwidth}
        \includegraphics[scale=0.31]{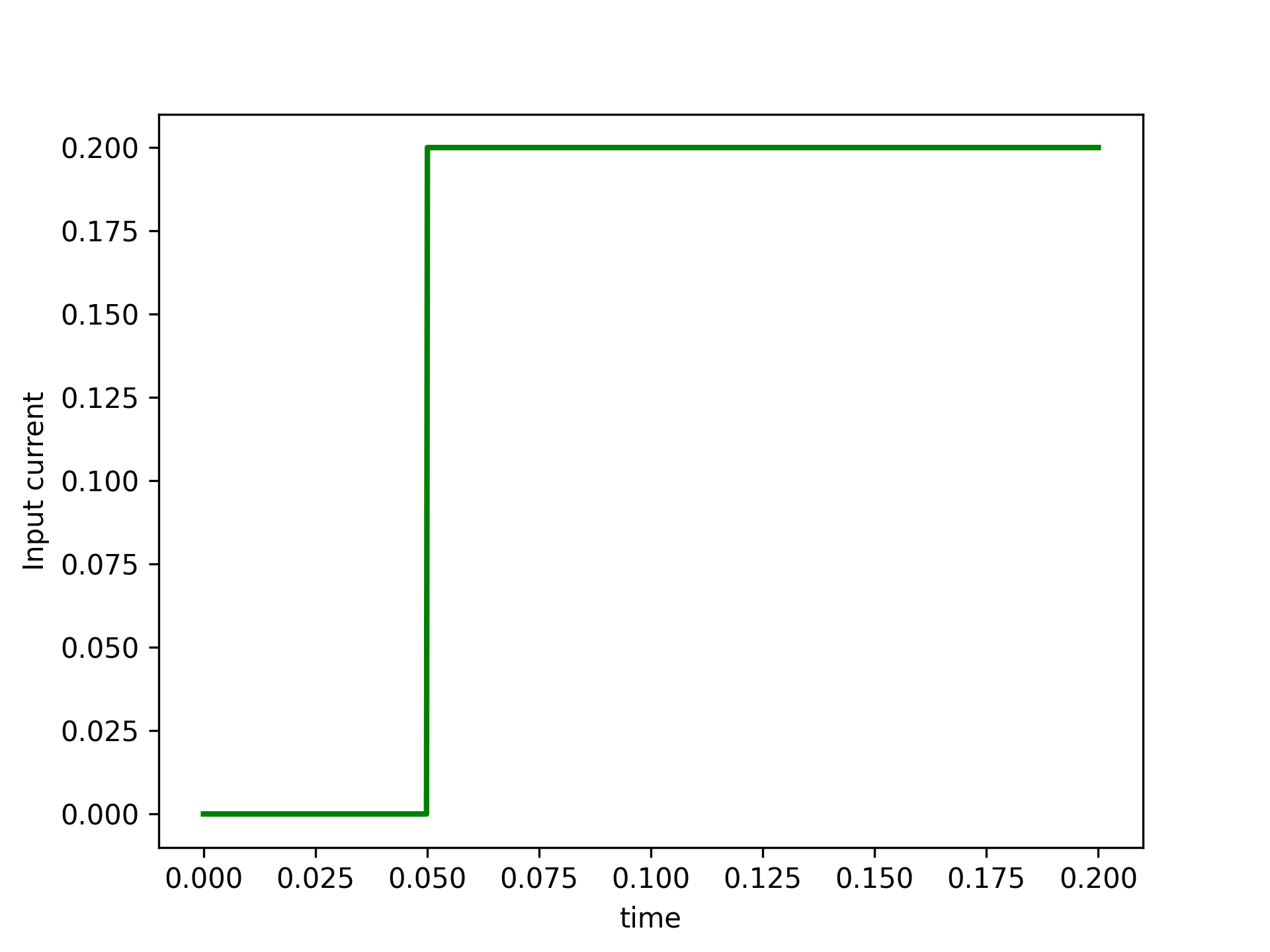}
        \caption{input current $I_{2}(t)$}
        \label{i2-3}
    \end{subfigure}
    \qquad
    \qquad
   \begin{subfigure}[b]{0.24\textwidth}
        \includegraphics[scale=0.30]{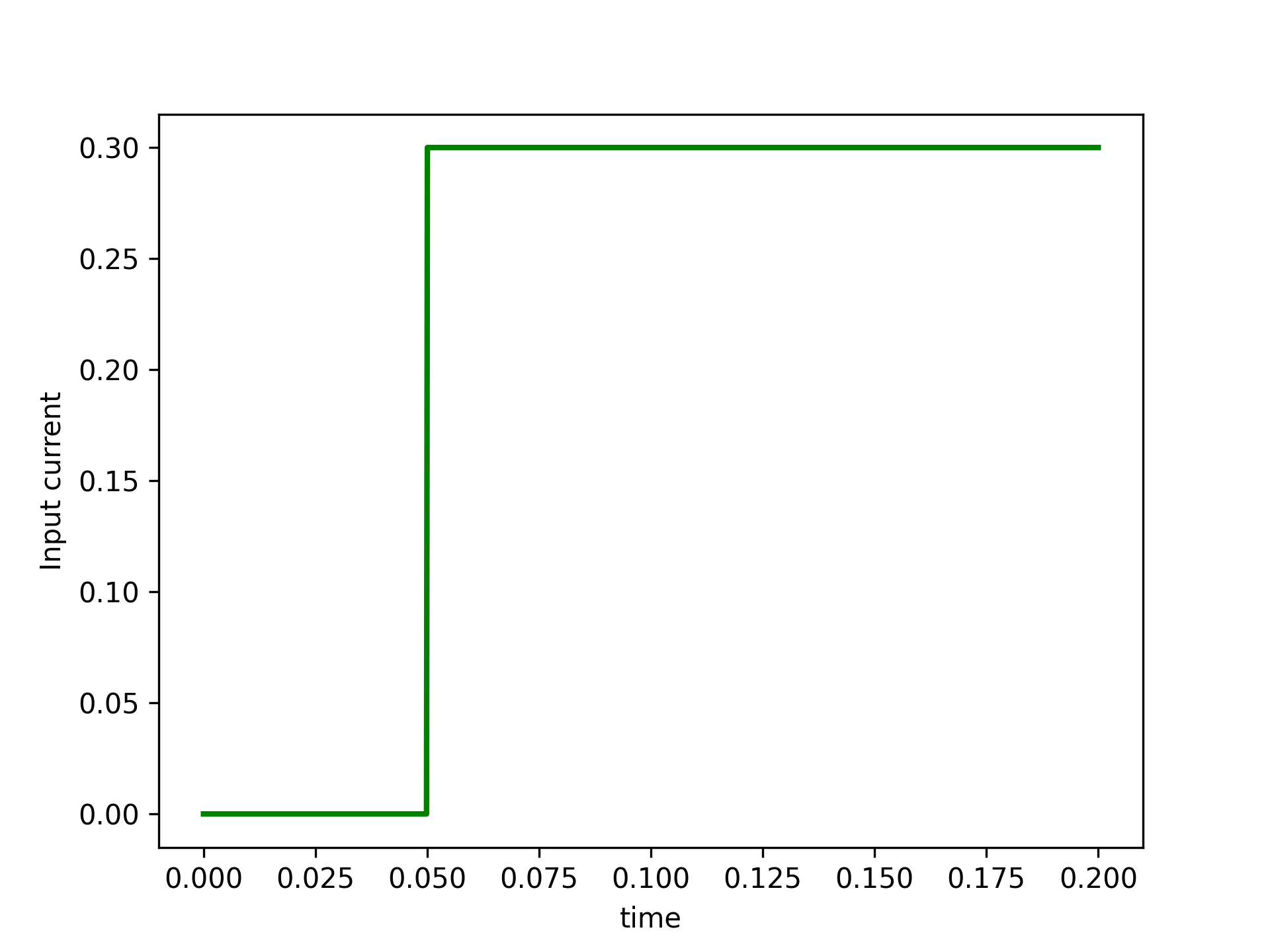}
        \caption{input current $I_{3}(t)$}
        \label{i2-4}
    \end{subfigure}
    \qquad
     \qquad
   \begin{subfigure}[b]{0.24\textwidth}
        \includegraphics[scale=0.30]{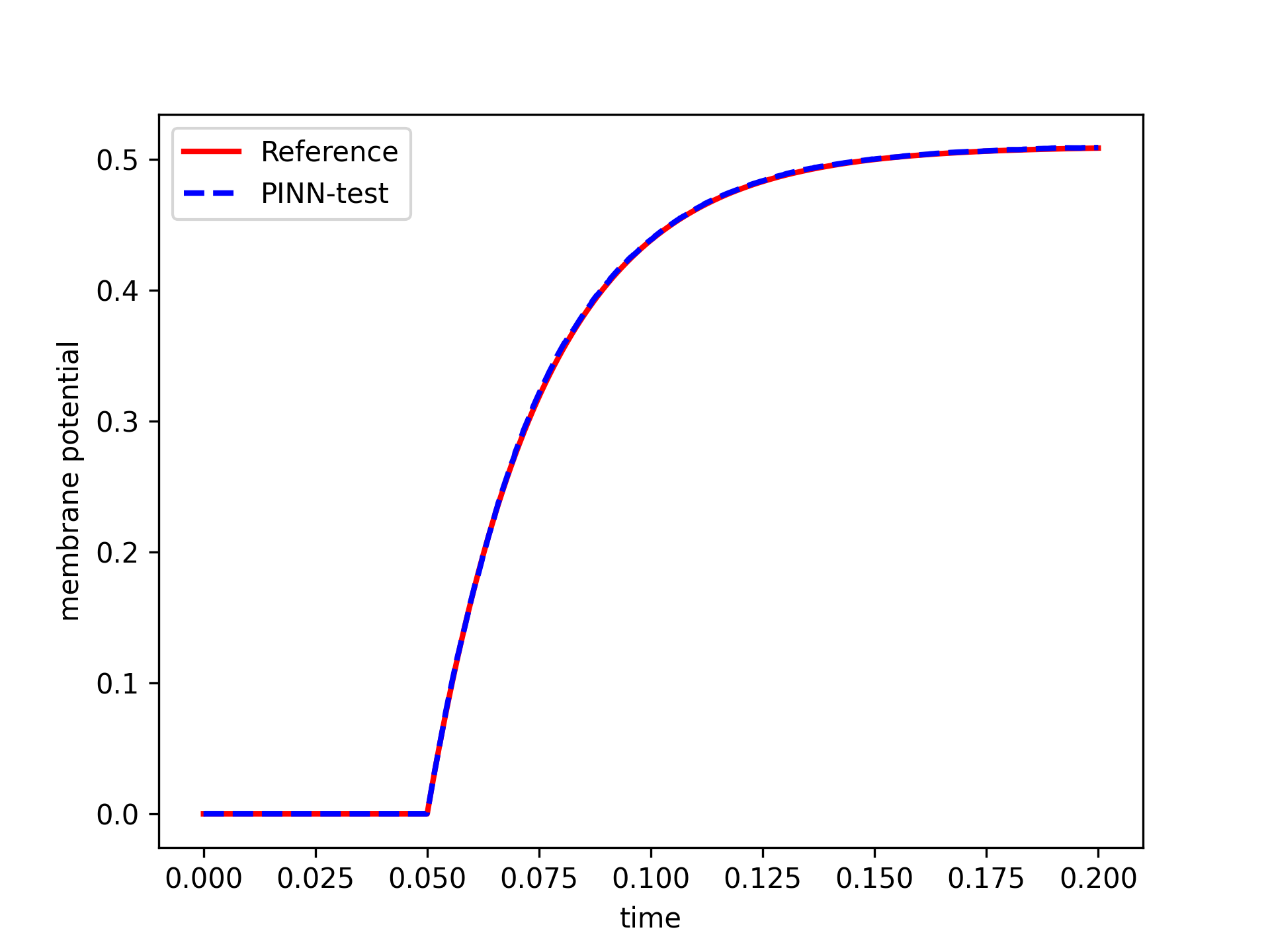}
        \caption{voltage response $V_{1}$}
        \label{i1-05}
    \end{subfigure}
    \qquad
    \qquad
    \begin{subfigure}[b]{0.26\textwidth}
        \includegraphics[scale=0.32]{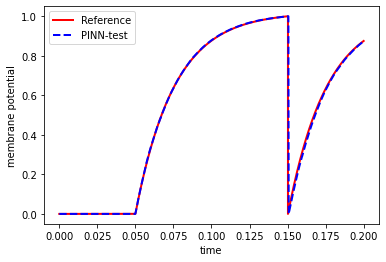}
        \caption{voltage response $V_{2}$}
        \label{i2-5}
    \end{subfigure}
    \qquad
    \qquad
    \begin{subfigure}[b]{0.24\textwidth}
        \includegraphics[scale=0.30]{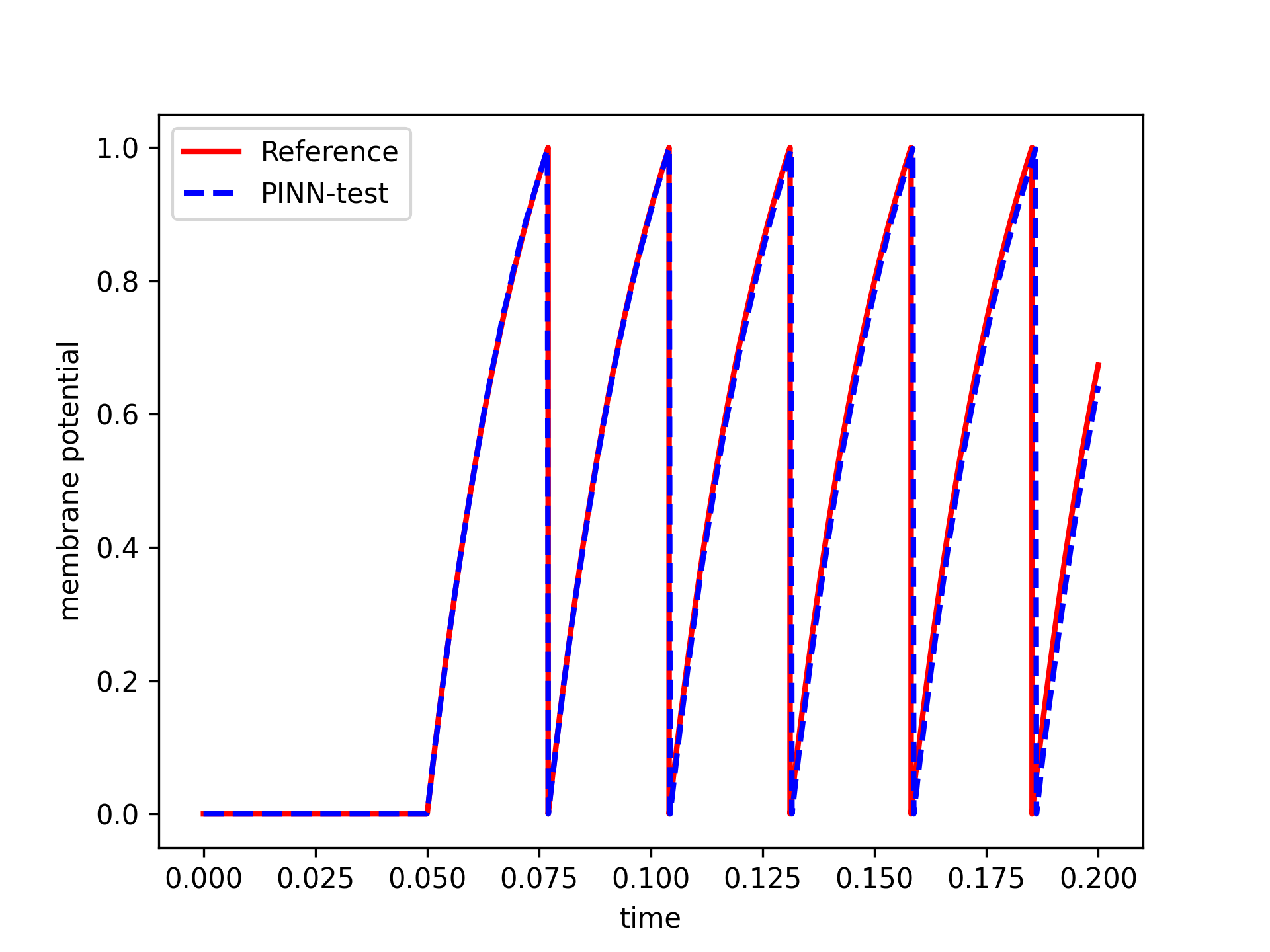}
        \caption{voltage response $V_{3}$}
        \label{i2-7}
    \end{subfigure}
    \qquad
    \qquad
    \begin{subfigure}[b]{0.27\textwidth}
        \includegraphics[scale=0.42]{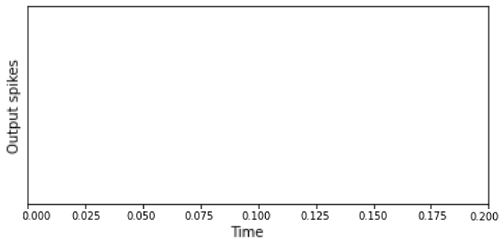}
        \caption{no output}
        \label{i2-8}
    \end{subfigure}
    \qquad
    \qquad
    \begin{subfigure}[b]{0.27\textwidth}
        \includegraphics[scale=0.42]{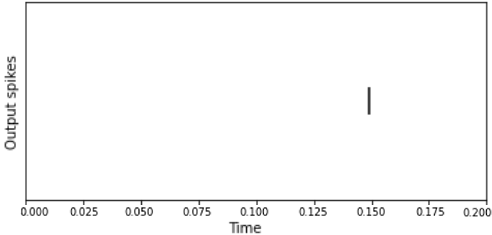}
        \caption{output of neuron}
        \label{i2-9}
    \end{subfigure}
    \qquad
     \qquad
    \begin{subfigure}[b]{0.21\textwidth}
        \includegraphics[scale=0.42]{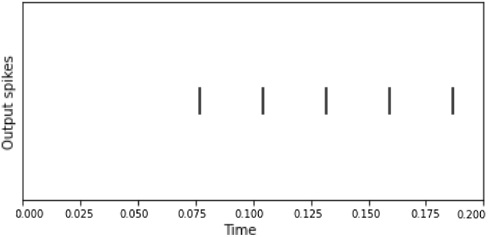}
        \caption{output of neuron}
        \label{i2-10}
    \end{subfigure}  
        \qquad
        \caption{LIF model; (a)-(c): input currents, (d)-(f): comparison of the PINN results with the reference solutions, and (g)-(i): output spikes.
        }
    \label{Fig44}
        \end{figure}

The approximate solutions of this model, when the threshold voltage is applied, have been shown in Figure \ref{Fig44}. In this case, for solving the problem in $[0, 0.2s]$, we have used $1000$ sub-intervals and $20$ points in each sub-interval. Furthermore, for solving the problem in each sub-interval, PINN has been used with the given architecture in Table \ref{tab}.

% {\color{red} Please check the third row of Figure \ref{Fig44}}
      
%\newpage
%\textbf{Appendix B: Reference solutions}\\

\section{Reference solutions}

We present a highly efficient and accurate spectral collocation method for obtaining the reference solutions of the fractional differential equations studied in the paper.
Consider the following system of FODEs
\begin{equation}\label{zeng-eq-1}
{}_CD_{0,t}^{\alpha_k} U_k(t)= f_k(t,U),\quad 0<t\le T,\quad k=1,2...,d,d\in \mathbb{N},
\end{equation}
subject to the initial condition $U_k(0)=\phi_k,f_k(t,U),U=(U_1,U_2,...,U_d)$.

We apply the fast time-stepping spectral collocation method to solve \eqref{zeng-eq-1}.
Divide the interval into subintervals $[t^{k-1},t^{k}]$ such that
$[0,T] = [t^0,t^1]\cup [t^1,t^2]\cup \cdot \cup [t^{N-1},t^{N}]$ and
$t^{N}=T$, $N$ is a positive number.
Let $p$ be a positive number and $\{x_j\}_{j=0}^p$ be the Legendre--Gauss--Lobatto (LGL) quadrature points
on the standard interval $[-1,1]$. Then, the LGL quadrature points $\{t^{k}_{j}\}_{j=0}^p$
in the interval $[t^{k-1},t^{k}]$ can be written as
$$t_j^k=\frac{1}{2}\left[(t^k-t^{k-1})x_j + t^{k-1}+t^k\right].$$
Denote $\Pi^p$ the interpolation operator as
$$\Pi^p u(t)|_{t\in [t^{k-1},t^k]}= \sum_{j=0}^pu(t^k_{j})\ell^k_j(t),\qquad
\ell^k_j(t)=\prod_{q=0,q\ne j}^{p}\frac{t-t^k_q}{t^k_j-t^k_q}.$$

Next, we show  the fast  calculation of ${}_CD_{0,t}^{\alpha} u(t)$,
the algorithm is illustrated by the following \textbf{Step 1)--Step 4)}.
%\textbf{Algorithm 1} Fast calculation of ${}_CD_{0,t}^{\alpha} u(t)$.
\begin{itemize}
  \item
\textbf{Step 1)} Find a sum-of-exponentials to approximate the kernel function
$k_{-\alpha}(t)$ as
\begin{equation}\label{soe}
k_{-\alpha}(t)=\frac{t^{-\alpha-1}}{\Gamma(-\alpha)}=\sum_{j=1}^{Q}w_je^{-\lambda_j t}
+ O(\varepsilon t^{-\alpha-1}),\quad t \in [t_1,T],
\end{equation}
where $\varepsilon$ is the  relative error, $Q>1$ is a positive integer,  and
\begin{equation*}\label{soe-err-wj}
\lambda_j=e^{y_j},\quad w_j=-\frac{\sin(\alpha \pi)}{\pi} \Delta y e^{(1+\alpha)y_j}.
\end{equation*}
Here, $y_j=y_{\min} + j\Delta y$, $\Delta y=(y_{\max}-y_{\min})/(Q-1)$,  and
\begin{equation*}\label{eq:c23-1}
  \begin{aligned}
& y_{\min}= (1+\alpha)^{-1}\ln(\varepsilon_{0}) - \ln(T), \quad \varepsilon_{0}=10^{-16},\\
& y_{\max}= \log\left(\frac{-\ln(\varepsilon_{0}) + (1+\alpha)\ln(t_1)}{0.5t_1}\right), \quad \varepsilon_{0}=10^{-16}.
  \end{aligned}
\end{equation*}

 \item
\textbf{Step 2)}
For $t>t^{n-1}$, divide
$${}_C D_{0,t}^{\alpha} u(t)= P.V.\int_{0}^{t} k_{-\alpha}(t-s) ( u(s)-u(0)) \ dx[s]$$
into two parts as
\begin{equation*}
\label{LH}
\begin{aligned}
{}_C D_{0,t}^{\alpha} u(t)
=  \underbrace{P.V.\int_{t^{n-1}}^{t} k_{-\alpha}(t-s)  u(s)\ dx[s]
-\frac{t^{-\alpha}u(0)}{\Gamma(1-\alpha)}}_{L^{\alpha} u}
+  \underbrace{\int_{0}^{t^{n-1}}k_{-\alpha}(t-s)  u(s)\ dx[s]}_{H^{\alpha} u},
\end{aligned}\end{equation*}
where $P.V.$ stands for the principal value with
$$P.V.\int_{a}^{t} k_{-\alpha}(t-s)u(s)\ dx[s]
= \frac{1}{\Gamma(1-\alpha)}\frac{\ dx[]}{\ dx[t]}\int_a^t(t-s)^{-\alpha}u(s)\ dx[s].$$

For $j=1,2,...,N$,
the local part $L u (t_j^n)$ is calculated by
$$ L^{\alpha} u (t_j^n)\approx L^{\alpha}u_j^n
= P.V.\int_{t^{n-1}}^{t_j^n} k_{-\alpha}(t_j^n-s)  (\Pi^pu)(s)\ dx[s]
-\frac{u(0)}{\Gamma(1-\alpha)}(t_j^n)^{-\alpha}.$$
%The history part $H^{\alpha} u^n$ is approximated by
%${}_FH^{\alpha} u^n$ in \textbf{Step 3)}.

 \item
\textbf{Step 3)}
By  \eqref{soe}, the history part $H^{\alpha} u(t_j^n),j=1,2,...,p$, is approximated by
\begin{equation}\label{FL1}\begin{aligned}
H^{\alpha} u(t_j^n)\approx &\int_{0}^{t^{n-1}}k_{-\alpha}(t_j^n-s) \Pi^p u(s)\ dx[s]\\
\approx &\sum_{k=1}^Qw_k \int_{0}^{t^{n-1}} e^{-\lambda_k(t_j^n-s)}\Pi^p u(s) \ dx[s]\\
= & \sum_{k=1}^Qw_ke^{-\lambda_k(t_j^n-t^{n-1})} \underbrace{\int_{0}^{t^{n-1}}
e^{-\lambda_k(t^{n-1}-s)}\Pi^p u(s) \ dx[s]}_{Y_k(t^{n-1})} \\
=& \sum_{k=1}^Qw_ke^{-\lambda_k(t_j^n-t^{n-1})} Y_k(t^{n-1}) = {}_FH^{\alpha} u_j^n,
\end{aligned}\end{equation}
where $Y_k(t)$ satisfies the following linear ODE
$$Y_k'(t)=-\lambda_kY_k(t) + \Pi^pu(t),\qquad Y_k(0) = 0,$$
which can be exactly solved by the following recurrence formula
$$Y_{k}(t^{n-1})=e^{-\lambda_{k}(t^{n-1}-t^{n-2})}Y_{j}(t^{n-2})+
\int_{t^{n-2}}^{t^{n-1}}e^{-\lambda_k(t^{n-1}-t)}\Pi^pu(t)\ dx[t],\quad Y_{j}(0)=0.$$

 \item  \textbf{Step 4)} Output $ {}_F\delta_t^{\alpha}  u^n_j = L^{\alpha} u^n_j + {}_FH^{\alpha} u^n_j$,
 $j=1,2,...,p$.
\end{itemize}

We are now in a position to develop the fast time-stepping collocation method for solving \eqref{zeng-eq-1}.

The fast time-stepping collocation method for \eqref{zeng-eq-1} is:
Given $U_{\ell,j}^k,k\le n-1,0\le j \le p,n\ge 1$,
find $U^n_{\ell,j}$, such that
\begin{equation}\label{zeng-eq-2}
{}_F\delta_t^{\alpha}  U^n_j = f(t^n_j,U^n_{\ell,j}),\qquad 1\le j \le p,\quad 1\le \ell \le d.
\end{equation}
The above equation can be rewritten as
\begin{equation}\label{zeng-eq-3}\begin{aligned}
P.V.\int_{t^{n-1}}^{t_j^n} k_{-\alpha_{\ell}}(t_j^n-s)(\Pi^pU_{\ell})(s)\ dx[s]-f_{\ell}(t^n,U^n_{\ell,j})\\
 =& -{}_FH^{\alpha_{\ell}} U^n_{\ell,j}
+\frac{(t_j^n)^{-\alpha_{\ell}}}{\Gamma(1-\alpha_{\ell})}\phi_{\ell},\quad 1\le j \le N,1\le\ell\le d.
\end{aligned}\end{equation}
The unknowns $u_j^n(1\le j \le N)$ in the above nonlinear system can be derived by the Newton iteration method.

We can rewrite the FOHH model \eqref{eq13-B-2}--\eqref{eq14-B}  can be solved by \eqref{zeng-eq-3}.

Some computational issues are addressed below.
\begin{itemize}
  \item[(i)] Generally speaking, the solution of the FODE, e.g., \eqref{zeng-eq-1},
  has a singularity at the origin ($t=0$ for  \eqref{zeng-eq-1}), we need to
  refine the mesh near the origin such that high accurate numerical solutions are obtained. For
  the numerical method \eqref{zeng-eq-3}, for $t\in [0,\Delta T]$, the mesh $t^n$ can be defined by   $t^n=\Delta T (n/M_1)^r$, $r\ge 1$. For $t\in [\Delta T,T]$, the
  uniform mesh is used which can be given by $t^n=t^{M_1} + (n-M_1)\Delta t$, $\Delta t=(T-\Delta T)/M_2$,
  $M_2=\lfloor(T-\Delta T)/(t^{M_1}-t^{M_1-1})\rfloor$.
  \item[(ii)] For the fractional LIF model with the current $I$ being piecewise constant,
  the solution also has a singularity at the points of the discontinuity of $I(t)$.
  We can use the strategy of (i) to refine the mesh near the discontinuity points when
  we use the method \eqref{zeng-eq-3} to solve the fractional LIF.

 % \item [(iii)]
\end{itemize}

\end{document}